\documentclass[journal]{IEEEtran}
\ifCLASSINFOpdf
\else
\fi
\usepackage{nomencl}
\usepackage{amssymb}
\usepackage{float}
\usepackage{graphicx}
\usepackage{subfigure}
\usepackage{etoolbox}
\usepackage{cite}
\usepackage{amsmath,amssymb,amsfonts}
\usepackage{amsmath,textcomp}
\usepackage{algorithmic}
\usepackage{textcomp}
\usepackage{xcolor}
\usepackage{multirow}
\usepackage{makecell}
\usepackage{booktabs}
\usepackage{dsfont}
\usepackage{textcomp}
\usepackage{threeparttable}
\usepackage{threeparttable}
\usepackage{upgreek}
\usepackage{mathrsfs}
\usepackage{CJKutf8}
\usepackage{xcolor}
\makenomenclature
\begin{document}
\title{Data-Driven Distributionally Robust Planning of Coastal Distribution Networks and PV-Storage-EV Stations for Low-Carbon Energy Transition}\author{Wenhao~Gao,\IEEEmembership{}
        Yongheng~Wang,\IEEEmembership{}
        Yangfan~Luo,\IEEEmembership{}
        Xinwei~Shen*\IEEEmembership{}
\thanks{Wenhao Gao, Yangfan Luo and Xinwei Shen are with Tsinghua Shenzhen International Graduate School, Tsinghua University.   (Corresponding: Xinwei Shen, sxw.tbsi@sz.tsinghua.edu.cn)}
\thanks{Yongheng Wang is with the Department of Electrical and Computer Engineering, The University of Hong Kong.}
\thanks{This work is supported in part by National Natural Science Foundation of China (52477102), Guangdong Basic and Applied Basic Research Foundation (2023A1515240055), Excellent Youth Basic Research Fund of Shenzhen (RCYX20231211090430053).}
}
\maketitle
\begin{abstract}
Coastal distribution networks (DNs) are becoming key infrastructure for the low-carbon energy transition as electric vehicle (EV) charging, photovoltaic (PV) generation, energy storage, and marine renewable procurement are increasingly coupled at the planning stage. Existing planning models, however, usually treat network expansion, PV-storage-EV station (PSES) deployment, low-carbon operation, and coastal asset degradation separately. This paper proposes a data-driven distributionally robust planning framework for coastal DNs and PSESs. The framework jointly optimizes radial network topology and PSES siting while incorporating salt-spray degradation costs, carbon-emission-flow (CEF)-based nodal carbon accounting, and tidal-generation-driven procurement costs. To represent heterogeneous uncertainty, the model combines continuous box uncertainty for load and PV availability with probability ambiguity over EV-charging scenarios generated from historical charging profiles. The resulting two-stage distributionally robust optimization problem is reformulated through exact relaxation of ESS mode binaries and KKT-based decomposition, and is solved by an inexact column-and-constraint generation (i-C\&CG) algorithm. {Case studies on a realistic 47-node coastal feeder show that salt-spray exposure can alter the selected topology and PSES siting, while the proposed planning scheme improves voltage quality, reduces annualized total cost, and reveals the cost--carbon trade-off of tidal procurement. The tailored i-C\&CG algorithm obtains the same optimal design as conventional C\&CG with about an order-of-magnitude reduction in solution time.}
\end{abstract}
\vspace{-0.3em}
\begin{IEEEkeywords}
Coastal distribution network, PV-storage-EV station, low-carbon energy transition, data-driven distributionally robust planning, salt-spray corrosion, carbon emission flow.
\end{IEEEkeywords}
\vspace{-1.5em}
\section*{Nomenclature}
\setlength{\IEEEdlabelindent}{0pt}
\vspace{-0.5em}

\noindent\textbf{Indices and Sets}
\vspace{-0.3em}
\begin{IEEEdescription}[\IEEEusemathlabelsep\IEEEsetlabelwidth{$P_t^{\mathrm{TG,max}},P_t^{\mathrm{TC,max}}$}]
\item[$i,j,k$] Index of nodes.
\item[$ij,jk$] Index of branches.
\item[$t$] Index of time intervals.
\item[$u$] Index of EVs.
\item[$a$] Index of areas.
\item[$s$] Index of EV charging scenarios.
\item[$\mathcal{I}$] Set of nodes.
\item[$\mathcal{L}$] Set of distribution lines.
\item[$\mathcal{T}$] Set of time intervals.
\item[$\mathcal{S}$] Set of EV charging scenarios.
\item[$\mathcal{U}_s$] Set of EV sessions.
\item[$\mathcal{A}$] Set of areas.
{\item[$\mathcal{I}^{\mathrm{Sub}},\mathcal{I}_a$] Set of substation nodes and nodes in area $a$.}
\end{IEEEdescription}

\vspace{-0.3em}
\noindent\textbf{Parameters}
\vspace{-0.3em}
\begin{IEEEdescription}[\IEEEusemathlabelsep\IEEEsetlabelwidth{$P_t^{\mathrm{TG,max}},P_t^{\mathrm{TC,max}}$}]
\item[$c_{ij}$] Construction cost of line $ij$.
\item[$c_i^{\mathrm{PSES}}$] Construction cost of PSES at node $i$.
\item[$\kappa_{\mathrm{Line}}^{\mathrm{Salt}},\kappa_{\mathrm{PSES}}^{\mathrm{Salt}}$] Salt-spray degradation coefficients of DN lines and PSESs.
\item[$c^{\mathrm{Ele}}$] Time-of-use electricity price.
\item[$D_i$] Virtual load at node $i$.
\item[$y_a^{\min},y_a^{\max}$] Minimum and maximum number of PSESs in area $a$.
\item[$R_{ij},X_{ij}$] Resistance and reactance of line $ij$.
\item[$S_{ij}^{\max}$] Maximum complex power of line $ij$.
{\item[$V^{\min},V^{\max}$] Minimum and maximum squared voltage amplitudes.}
\item[$P^{\mathrm{Sub,max}}$] Maximum active power of substation.
\item[$Q^{\mathrm{Sub,max}}$] Maximum reactive power of substation.
\item[$\bar{P}_{i,t}^{\mathrm{PV}}$] Maximum active power of PV at node $i$ in time interval $t$.
\item[$\mu_i^{\mathrm{ESS,ch}},\mu_i^{\mathrm{ESS,dch}}$] Charging and discharging efficiencies of ESS at node $i$.
\item[$P_i^{\mathrm{ESS,min}},P_i^{\mathrm{ESS,max}}$] Minimum and maximum active power of ESS at node $i$.
\item[$E_i^{\mathrm{ESS,min}},E_i^{\mathrm{ESS,max}}$] Minimum and maximum energy capacity of ESS at node $i$.
\item[$\tau_u$] Set of time intervals of EV charging.
\item[$E_{u,0}^{\mathrm{EV}},E_u^{\mathrm{EV}}$] Initial and target energy of the $u$-th EV.
\item[$E_u^{\mathrm{EV,min}},E_u^{\mathrm{EV,max}}$] Minimum and maximum energy capacity of the $u$-th EV.
\item[$P_u^{\mathrm{EV,max}}$] Maximum charging power of the $u$-th EV.
\item[$\bar{P}_{i,t}^{\mathrm{Grid}},\bar{Q}_{i,t}^{\mathrm{Grid}}$] Reference active and reactive grid power.
\item[$\bar{P}_{i,t}^{\mathrm{PV},0}$] Reference active power of PV.
\item[$\Delta P_{i,t}^{-},\Delta P_{i,t}^{+}$] Minimum and maximum fluctuation ranges of active grid power.
\item[$\Delta Q_{i,t}^{-},\Delta Q_{i,t}^{+}$] Minimum and maximum fluctuation ranges of reactive grid power.
\item[$\Delta P_{i,t}^{\mathrm{PV},-},\Delta P_{i,t}^{\mathrm{PV},+}$] Minimum and maximum fluctuation ranges of PV output.
\item[$c_a^{\mathrm{CEF}}$] Carbon tax of area $a$.
\item[$c_t^{\mathrm{TG}},c_t^{\mathrm{TC}}$] Electricity prices from thermal and tidal-current units in time interval $t$.
\item[$P_t^{\mathrm{TG,max}},P_t^{\mathrm{TC,max}}$] Maximum outputs of thermal and tidal-current units in time interval $t$.
{\item[$\Delta t$] Duration of each time interval.}
\end{IEEEdescription}

\vspace{-0.3em}
\noindent\textbf{Variables}
\vspace{-0.3em}
\begin{IEEEdescription}[\IEEEusemathlabelsep\IEEEsetlabelwidth{$P_t^{\mathrm{TG,max}},P_t^{\mathrm{TC,max}}$}]
\item[$y_i,z_{ij}$] Binary variables associated with PSES at node $i$ and line $ij$.
\item[$P_{ij,t,s},Q_{ij,t,s}$] Active and reactive power flow in line $ij$ in time interval $t$ for scenario $s$.
\item[$e_{i,t}$] Carbon intensity at node $i$ in time interval $t$.
\item[$P_{i,t}^{\mathrm{Load}}$] Active load at node $i$ in time interval $t$.
\item[$P_t^{\mathrm{TG}},P_t^{\mathrm{TC}}$] Active power generated by thermal and tidal-current units.
\item[$y_{ij}$] Binary variable associated with parent-child relation between nodes $i$ and $j$.
\item[$F_{ij}$] Virtual power flow in line $ij$.
\item[$P_{i,t}^{\mathrm{Grid}},Q_{i,t}^{\mathrm{Grid}}$] Active and reactive load power at node $i$ in time interval $t$.
\item[$P_{i,t,s}^{\mathrm{EV}}$] EV charging power at node $i$ in time interval $t$ for scenario $s$.
\item[$P_{i,t,s}^{\mathrm{PV}}$] Active PV output at node $i$ in time interval $t$ for scenario $s$.
\item[$P_{i,t,s}^{\mathrm{ESS,ch}},P_{i,t,s}^{\mathrm{ESS,dch}}$] Charging and discharging power of ESS at node $i$, time $t$, and scenario $s$.
{\item[$V_{i,t,s}$] Squared voltage magnitude at node $i$ in time interval $t$ for scenario $s$.}
\item[$P_{i,t,s}^{\mathrm{Sub}},Q_{i,t,s}^{\mathrm{Sub}}$] Active and reactive power output of substation at time $t$ and scenario $s$.
\item[$E_{i,t,s}^{\mathrm{ESS}}$] Energy stored in ESS at node $i$ in time interval $t$ for scenario $s$.
\item[$\gamma_{i,t,s}^{\mathrm{ESS,ch}},\gamma_{i,t,s}^{\mathrm{ESS,dch}}$] Binary variables for ESS charging and discharging states.
\item[$P_{u,t,s}^{\mathrm{EV,sch}}$] Scheduled charging power of the $u$-th EV in time interval $t$ for scenario $s$.
\item[$P_{i,t,s}^{\mathrm{EV,sch}},P_{i,t,s}^{\mathrm{EV,ran}}$] Aggregated scheduled and random EV charging power at node $i$.
\item[$e_{i,t}^{\mathrm{G}}$] Carbon emission intensity of substation node in time interval $t$.
\end{IEEEdescription}

\vspace{-0.3em}
\noindent\textbf{Uncertainty and Distributions}
\vspace{-0.3em}
\begin{IEEEdescription}[\IEEEusemathlabelsep\IEEEsetlabelwidth{$P_t^{\mathrm{TG,max}},P_t^{\mathrm{TC,max}}$}]
\item[$\xi$] Continuous uncertainty realization.
\item[$\Xi$] Box-type uncertainty set for $\xi$.
\item[$\zeta_s^{\mathrm{EV}}$] EV charging profile for scenario $s$.
\item[$\pi$] Probability vector over EV scenarios.
\item[$\pi_s^0$] Reference probability for EV charging scenario $s$.
\item[$\Pi^{\mathrm{EV}}$] Ambiguity set of $\pi$ with combined $1$-norm and $\infty$-norm constraints.
\item[$\theta_1,\theta_\infty$] Radii of the $1$-norm and $\infty$-norm ambiguity constraints.
\item[$\alpha_1,\alpha_\infty$] Confidence levels used to set $\theta_1$ and $\theta_\infty$.
\end{IEEEdescription}

\vspace{-0.5em}
\section{Introduction}

\IEEEPARstart{C}{oastal} distribution networks are becoming important carriers of the low-carbon energy transition.
In coastal cities, distribution networks (DNs) are increasingly coupled with electric vehicle (EV) charging demand, distributed photovoltaic (PV) generation, energy storage systems (ESSs), and marine renewable resources such as tidal currents \cite{lee2014electric}.
This coupling changes the planning task.
DN topology affects the delivery of charging power and local renewable generation.
PV and ESS shape voltage profiles and losses \cite{wang2026two}.
Tidal procurement changes the carbon intensity and price of grid-side energy.
Therefore, network expansion, charging infrastructure deployment, storage operation, renewable procurement, and carbon-emission management need to be coordinated at the planning stage rather than optimized as separate decisions \cite{li2022techno}.

Existing studies provide useful foundations, but they do not fully address this coastal DN planning problem under the low-carbon transition.
Coordinated planning of DNs, EV charging facilities, and distributed generation resources has been widely studied \cite{zhang2023restoration,chen2024coordinated,zhang2024bilevel,shi2025fully}.
Active DN planning with PV and ESS flexibility has also been used to improve system performance under increasing electrification \cite{wang2023joint}.
At the station level, PV--storage--EV stations (PSESs) can coordinate local PV generation, storage dispatch, and EV charging to improve economic and emission performance \cite{sun2022economic,wang2025tri,liu2024transforming}.
However, most PSES planning studies treat station siting separately from DN topology investment, while many DN planning studies simplify charging stations as loads and do not model the internal PV-storage flexibility of integrated stations \cite{wang2026integrated,ren2018coordinated}.

A second limitation is that the coastal low-carbon context is often simplified.
{Existing coastal studies have addressed generation--storage assessment and sizing in communities and microgrids \cite{li2022techno,bakeer2023optimal}, as well as robust and stochastic planning with marine renewables and desalination \cite{wang2023robust,duah2025marine}.
Alongside these developments, salt-spray degradation introduces an additional cost consideration for DN--PSES planning.}
Unlike inland feeders, coastal DNs are exposed to chloride-rich humidity and salt spray, which accelerate the degradation of metallic lines and station components and increase life-cycle costs \cite{ISO9223-2012,yang2024review,alcantara2015airborne,gao2026life}.
Meanwhile, coastal cities may procure electricity from marine renewable resources, whose output and price affect both operating cost and carbon intensity.
Low-carbon planning also requires more than aggregate emission calculation, because nodal carbon intensity depends on power-flow patterns and supply mix.
Carbon emission flow (CEF) offers a power-flow-consistent way to allocate carbon intensity across network nodes and loads \cite{cheng2018planning,li2024electricity}.
{These features are closely related, motivating their joint consideration in DN--PSES planning.}

\begin{table*}[!t]
\centering
\begingroup

\setlength{\tabcolsep}{3pt}
\renewcommand{\arraystretch}{1.16}
\caption{Comparison of research scope, uncertainty modeling, and solution methods}
\label{tab:coastal_literature}
\begin{tabular}{@{}
>{\raggedright\arraybackslash}p{0.08\textwidth}
>{\raggedright\arraybackslash}p{0.205\textwidth}
>{\raggedright\arraybackslash}p{0.17\textwidth}
>{\raggedright\arraybackslash}p{0.105\textwidth}
>{\raggedright\arraybackslash}p{0.23\textwidth}
>{\raggedright\arraybackslash}p{0.15\textwidth}@{}}
\hline\hline
\multirow{2}{*}{Refs.}
& \multicolumn{2}{c}{Research scope}
& \multicolumn{2}{c}{Uncertainty modeling}
& \multirow{2}{*}{\makecell[l]{Solution method}} \\
\cline{2-5}
& Research focus
& Coastal/marine factors
& Method
& Representation
& \\ \hline

\cite{li2022techno}
& Coastal community:\newline generation--storage assessment
& Offshore wind; tidal energy
& Deterministic
& Given time-series inputs
& Time-series simulation
\\[4pt]

\cite{bakeer2023optimal}
& Coastal microgrid:\newline generation--storage sizing
& Coastal wind and solar 
profiles
& Deterministic
& Given time-series inputs
& Metaheuristic
\\[4pt]

\cite{wang2023robust}
& Offshore microgrid:\newline generation--storage planning
& Tidal generation; desalination
& RO
& Budgeted sets (load, tidal height)
& C\&CG
\\[4pt]

\cite{duah2025marine}
& Coastal community:\newline multi-energy capacity planning
& Offshore wind; wave energy; desalination
& SP
& Fixed-probability demand and resource scenarios
& Direct solve
\\ \hline

\cite{xiang2024robust,gao2024distributionally}
& DN planning with charging stations or multi-energy stations
& Not explicitly modeled
& DRO
& Combined $1$/$\infty$-norm \cite{xiang2024robust};\newline
  KL divergence \cite{gao2024distributionally}
& C\&CG
\\[4pt]

\cite{nguyen2022decentralized}
& DN--EVCS operational coordination
& Not explicitly modeled
& DRO
& Wasserstein ambiguity set
& ADMM
\\ \hline

Proposed
& Coastal DN--PSES planning:\newline joint topology and siting
& Salt-spray costs; tidal procurement
& DRO
& EV: combined $1$/$\infty$-norm; \newline
Load/PV: box uncertainty sets
& KKT reformulation; tailored i-C\&CG
\\
\hline\hline
\end{tabular}
\endgroup
\end{table*}

{Table~\ref{tab:coastal_literature} compares representative coastal studies and related DN applications in terms of research scope, uncertainty modeling, and solution methods.
Building on this literature, the main research gap is the lack of a data-driven planning framework that jointly captures coupled elements of coastal DNs and PSESs for low-carbon energy transition: DN topology investment, PSES siting and operation, coastal stressor cost, and network-resolved carbon and procurement effects.
This integrated formulation is designed to avoid biases that may occur when the above factors are treated separately: underestimating coastal degradation costs, misrepresenting the carbon effect of siting decisions, or selecting charging-station locations that are feasible for EV demand but suboptimal for network operation and low-carbon performance.}

{Uncertainty further complicates this planning problem because the uncertain quantities have different data structures.}
EV charging demand is naturally represented by spatio-temporal scenarios learned from historical charging profiles, but the scenario probabilities may be misspecified under limited samples, nonstationary travel behavior, and price-responsive charging \cite{zhu2025real}.
By contrast, base load and PV availability are continuous uncertain quantities that can be described by tractable uncertainty intervals \cite{gao2022integrated}.
Stochastic programming (SP) depends on an accurate probability distribution, whereas robust optimization (RO) may be too conservative for long-term investment decisions \cite{yang2024stochastic}.
Distributionally robust optimization (DRO) provides a middle ground by optimizing against the worst-case distribution within a data-calibrated ambiguity set \cite{nguyen2022decentralized,xiang2024robust,gao2024distributionally}.
This motivates a hybrid DRO formulation that uses historical EV data to build scenario-probability ambiguity while retaining box uncertainty for continuous load and PV deviations.

The resulting planning model is computationally challenging because topology binaries, PSES siting binaries, ESS operating modes, power-flow constraints, and uncertainty responses are coupled.
Standard column-and-constraint generation (C\&CG) is widely used for two-stage RO and DRO problems \cite{zeng2013solving}.
However, when the master problem contains both topology and siting decisions, solving it to high precision at every iteration can dominate the computational burden.
The inexact C\&CG (i-C\&CG) method can reduce this burden by allowing controlled inexact master solutions and using backtracking to recover valid bounds \cite{tsang2023inexact}.
Its role in coastal DN topology and PSES siting co-planning with hybrid uncertainty remains insufficiently studied.

\begin{figure}[htbp]
    \centering
    \includegraphics[width=0.47\textwidth]{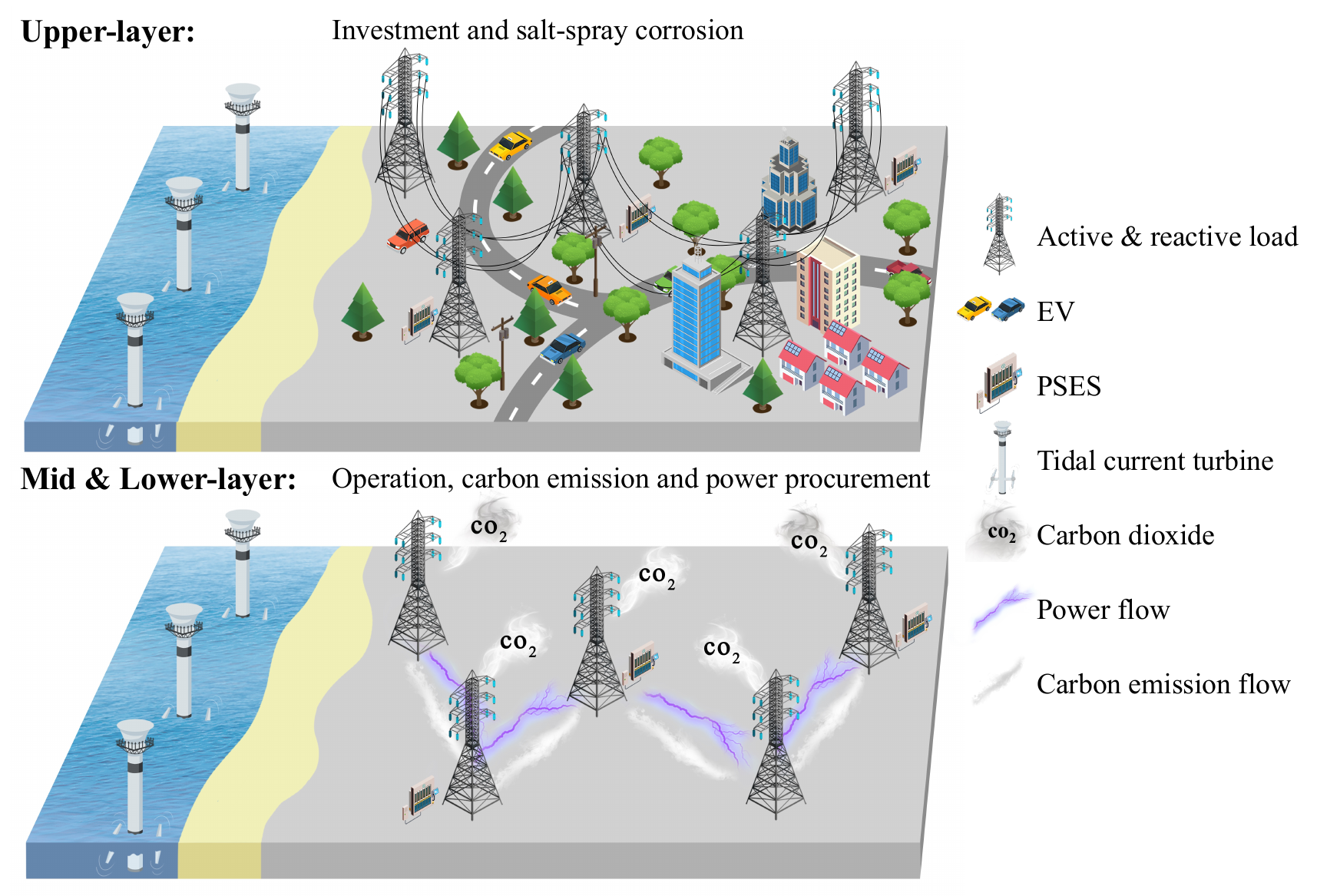}
    \vspace{-1.0em}
    \caption{{Schematic of the coastal DN–PSES planning context, showing network and station investment with salt-spray degradation costs (upper layer) and network operation, carbon emission flow accounting, and tidal energy procurement (middle and lower layers).}}\label{fig1}
\end{figure}

{This paper addresses the above gap by developing a data-driven distributionally robust planning framework for coastal DNs and PSESs, as illustrated in Fig.~\ref{fig1}.
The four contributions are summarized below, with the first three addressing technical developments and the fourth presenting case-based coastal planning insights.}

\textbf{\textit{1) Coastal DN--PSES planning model:}}
We formulate a coastal DN--PSES co-planning model that jointly determines radial DN topology and PSES siting, instead of treating charging stations as fixed loads or optimizing PSESs outside the network expansion problem. 
The model incorporates salt-spray degradation costs, tidal-generation-driven procurement costs, and CEF-based nodal carbon accounting, thereby connecting coastal physical stressors with low-carbon planning decisions.

\textbf{\textit{2) Hybrid data-driven DRO formulation:}}
We develop a two-stage DRO formulation that treats EV charging uncertainty through discrete historical-data-driven scenarios with ambiguous probabilities, while modeling load and PV deviations through continuous box uncertainty.
This formulation reflects the different uncertainty structures of electricity demand, renewable availability, and transport charging behavior.

{\textbf{\textit{3) Tractable reformulation and solution method:}}
We relax the ESS mode binaries without loss of optimality under the proposed formulation, which enables a convex mid-layer reformulation and KKT-based decomposition.
An i-C\&CG algorithm is then tailored to solve the resulting coupled DN topology and PSES siting problem efficiently.}

{\textbf{\textit{4) Coastal planning insights:}}
Comparative studies on a 47-node coastal DN identify differences in the selected topology and station locations with and without salt-spray degradation costs, and quantify the voltage and loss benefits of PSES deployment. Tidal-procurement sensitivity analyses further characterize the cost–carbon implications of marine renewable supply. These case-based insights highlight the relevance of coastal asset costs, local PV–storage flexibility, and procurement conditions when assessing DN–PSES plans.}

The remainder of this paper is organized as follows.
Section II presents the proposed model.
Section III develops the solution method based on i-C\&CG.
Section IV reports case studies and sensitivity analyses.
Section V concludes the paper.

\section{Model Formulation}
\subsection{Objective Function}
The objective function is designed for coastal DN--PSES co-planning under the low-carbon transition. It includes DN and PSES investment costs, the incremental cost induced by salt spray, DN power-loss cost, carbon-emission cost, and power-procurement cost. The overall objective function is given by:
\vspace{-0.5em}
\begin{align}
\min_{\mathbf{z},\mathbf{y}} \Big\{ \;&
\underbrace{
C^{\mathrm{Inv}}(\mathbf{z},\mathbf{y})+ C^{\mathrm{Salt}}(\mathbf{z},\mathbf{y})
}_{\substack{\text{Upper-layer: planning}\\[-0.2ex]\text{(DN topology \& PSES siting)}}}
\notag\\
&+\underbrace{
\max_{\xi \in \Xi}\;\sup_{\pi \in \Pi^{\mathrm{EV}}}\;
\mathbb{E}_{s\sim \pi}\!\left[
\min_{\mathbf{x}_s \in \mathcal{Y}_s(\mathbf{z},\mathbf{y};\xi,\zeta_s^{\mathrm{EV}})}
C^{\mathrm{Ope}}(\mathbf{x}_s)
\right]
}_{\substack{\text{Mid-layer: DRO operation}\\[-0.2ex]\text{(worst-case expected DN operating loss)}}}
\notag\\
&+\underbrace{
\min_{\mathbf{e},\,\mathbf{p}^{\prime}} \;
\Big(C^{\mathrm{Carb}}(\mathbf{e})+C^{\mathrm{Powr}}(\mathbf{p}^{\prime})\Big)
}_{\substack{\text{Lower-layer: carbon emission and procurement}\\[-0.2ex]\text{(decoupled and solved post-optimization)}}}
\Big\}.
\label{eq_Obj}
\end{align}

In \eqref{eq_Obj}, \(\mathbf{x}_s\) stacks the DN operating variables under EV-charging scenario \(s\) (line flows, nodal voltages, PV/ESS dispatch, and substation injections), and \(\mathcal{Y}_s(\mathbf{z},\mathbf{y};\xi,\zeta_s^{\mathrm{EV}})\) denotes the corresponding scenario-dependent feasible set of \(\mathbf{x}_s\) induced by the network operation constraints for given planning decisions \((\mathbf{z},\mathbf{y})\) and uncertainty realization \((\xi,\zeta_s^{\mathrm{EV}})\). The lower-layer decision variables are \(\mathbf{e}:=\{e_{i,t}\}\) and \(\mathbf{p}^{\prime}:=\{P_t^{\mathrm{TG}},P_t^{\mathrm{TC}}\}_{t\in \mathcal{T}}\), which are used only for carbon-emission accounting and power procurement.

\noindent \emph{1) Investment and Salt-spray Degradation Cost:}
\begin{flalign}
C^\mathrm{Inv}(\mathbf{z},\mathbf{y}) = \sum_{ij\in \mathcal{L}} c_{ij} z_{ij}+ \sum_{i\in \mathcal{I}} c_i^\mathrm{PSES} y_i,
\label{eq_Inv}
\end{flalign}
\vspace{-1.2em}
\begin{flalign}
C^\mathrm{Salt}(\mathbf{z},\mathbf{y}) = \sum_{ij\in \mathcal{L}} \kappa_\mathrm{Line}^\mathrm{Salt} c_{ij} z_{ij} + \sum_{i\in \mathcal{I}} \kappa_\mathrm{PSES}^\mathrm{Salt} c_i^\mathrm{PSES} y_i.
\label{eq_Salt}
\end{flalign}

The investment term \eqref{eq_Inv} represents the annualized cost of DN lines and PSESs, where the binary decision variables $z_{ij}$ and $y_i$ indicate whether line $ij$ and a PSES at node $i$ are constructed. The parameters $c_{ij}$ and $c_i^\mathrm{PSES}$ correspond to the associated annualized construction costs. Salt-spray corrosion is a distinctive stressor for coastal distribution systems, accelerating the degradation of exposed metallic components of DN lines and PSESs \cite{alcantara2015airborne}. To reflect its planning-level economic impact without introducing additional state variables, the annual corrosion-related surcharge is modeled as a linear fraction of the corresponding investment cost in \eqref{eq_Salt}, calibrated using ISO 9223:2012 \cite{ISO9223-2012}.
\begin{flalign*}
\kappa_\mathrm{Line}^\mathrm{Salt}=0.0175\,
\left(d^{\mathrm{Zn}}\right)^{-1}\left(S_d^{\mathrm{Cl}}\right)^{0.57} \exp (0.008 \mathrm{RH}+0.085 \vartheta),\end{flalign*}
\vspace{-2em}
\begin{flalign*}
\kappa_\mathrm{PSES}^\mathrm{Salt}=0.0018\,\left(d^{\mathrm{Al}}\right)^{-1}\left(S_d^{\mathrm{Cl}}\right)^{0.60} \exp (0.020 \mathrm{RH}+0.094 \vartheta),\end{flalign*}
where $\kappa_\mathrm{Line}^\mathrm{Salt}$ and $\kappa_\mathrm{PSES}^\mathrm{Salt}$ denote the degradation coefficients derived from chloride-driven dose--response functions using annual-average relative humidity $\mathrm{RH}$, ambient temperature $\vartheta$, chloride deposition rate $S_d^{\mathrm{Cl}}$, protective thicknesses $d^{\mathrm{Zn}}$ and $d^{\mathrm{Al}}$. Salt-sensitive exposure of DN lines is mainly associated with zinc-galvanized steel hardware rather than conductors, whereas PSES exposure is dominated by aluminum parts such as PV module frames \cite{IEEE_C13580_2012}. Physically, the coefficients represent corrosion penetration normalized by the allowable protective thickness: faster corrosion or thinner protection implies higher salt-spray-induced maintenance costs. For coastal districts, annual-average meteorology is assumed spatially homogeneous, so the coefficients are treated as constants, consistent with salinity-driven increases in zinc/aluminum corrosion rates in subtropical coastal environments \cite{panchenko2017prediction}.

\noindent \emph{2) Operation Cost:}
\begin{align}
{C^\mathrm{Ope}(\mathbf{x}_s) = \eta^\mathrm{Grid} \sum_{t \in \mathcal{T}} \sum_{ij \in \mathcal{L}} c^\mathrm{Ele}_t R_{ij}\left(P_{ij,t,s}^2+Q_{ij,t,s}^2\right) \Delta t.}
\label{eq_Ope}
\end{align}
\vspace{-1em}

The operation term \eqref{eq_Ope} represents the DN cost caused by power losses. Here, $\eta^\mathrm{Grid}$ denotes the annual number of operating days, which is typically assumed to span the entire year. $c^\mathrm{Ele}$ reflects the TOU electricity tariff variation between peak and off-peak periods. $P_{ij,t,s}$ and $Q_{ij,t,s}$ denote the active and reactive power flows on line $ij$ for scenario $s$, and $R_{ij}$ is the resistance of the corresponding line \cite{wang2023joint}.

\noindent \emph{3) Carbon Emission and Power Procurement Cost:}
\begin{align}
{C^\mathrm{Carb}(\mathbf{e}) = \eta^\mathrm{Grid} \sum_{t \in \mathcal{T}}\sum_{a \in \mathcal{A}} \sum_{i \in \mathcal{I}_a} c_{a}^\mathrm{CEF} e_{i,t} P_{i,t}^\mathrm{Load} \Delta t,}
\label{eq_Carb}
\end{align}
\vspace{-1em}
\begin{align}
{C^\mathrm{Powr}(\mathbf{p}^{\prime}) = \eta^\mathrm{Grid} \sum_{t \in \mathcal{T}}\left(c_t^\mathrm{TG} P_t^\mathrm{TG}+c_t^\mathrm{TC} P_t^\mathrm{TC}\right) \Delta t.}
\label{eq_Powr}
\end{align}

The cost terms \eqref{eq_Carb} and \eqref{eq_Powr} represent carbon emission and power procurement costs for the coastal DN, respectively. $c_{a}^\mathrm{CEF}$ denotes the carbon tax in area $a$, $e_{i,t}$ is the carbon intensity at node $i$ during each time interval, and $P_{i,t}^\mathrm{Load}$ represents the active load at the corresponding node. For power procurement, $P_t^\mathrm{TG}$ and $P_t^\mathrm{TC}$ indicate the active power generated by thermal units and tidal current units, while $c_t^\mathrm{TG}$ and $c_t^\mathrm{TC}$ denote the corresponding electricity prices for each time interval.

\vspace{-1.5em}
\subsection{Constraints}
\noindent (1) \textit{Radial Topology Constraints:}
\begin{flalign}
\sum_{ij \in \mathcal{L}} z_{ij}=|\mathcal{I}|-|\mathcal{I}^\mathrm{Sub}|,
\label{eq_topology1}
\end{flalign}
\vspace{-1.5em}
\begin{flalign}
y_{ij}+y_{ji}=z_{ij}, \quad \forall ij\in \mathcal{L},
\label{eq_topology2}
\end{flalign}
\vspace{-1.5em}
\begin{flalign}
y_{ij}=0, \quad \forall i\in \mathcal{I}^\mathrm{Sub},\forall ij\in \mathcal{L},
\label{eq_topology3}
\end{flalign}
\vspace{-1.5em}
\begin{flalign}
\sum_{ij\in \mathcal{L}} y_{ij}=1, \quad \forall i \in \mathcal{I} \backslash \mathcal{I}^\mathrm{Sub},
\label{eq_topology4}
\end{flalign}
\vspace{-1.5em}
\begin{flalign}
\sum_{ij\in \mathcal{L}} F_{ij}+D_i=\sum_{ki\in \mathcal{L}} F_{ki}, \quad \forall i\in \mathcal{I},
\label{eq_topology5}
\end{flalign}
\vspace{-1.5em}
\begin{flalign}
\left|F_{ij}\right| \leq z_{ij} M, \quad \forall ij \in \mathcal{L}.
\label{eq_topology6}
\end{flalign}

Cons. \eqref{eq_topology1}--\eqref{eq_topology6} define the DN radial topology constraints \cite{wang2020radiality}. Cons. \eqref{eq_topology1} ensures that the number of lines equals the number of non-root nodes. Cons. \eqref{eq_topology2}--\eqref{eq_topology4} impose spanning tree constraints. Cons. \eqref{eq_topology2} guarantees each connected line maintains a single parent-child relationship: if node $j(i)$ is the parent of node $i(j)$, then $y_{ij}$ ($y_{ji}$) equals 1. Cons. \eqref{eq_topology3} and \eqref{eq_topology4} enforce that root nodes have no parent node, whereas all other nodes are associated with exactly one parent.
Cons. \eqref{eq_topology5}--\eqref{eq_topology6} incorporate single-commodity flow constraints to prevent pseudo-root nodes. $F_{ij}$ is the fictitious flow and $D_i$ is the fictitious demand. Cons. \eqref{eq_topology5} ensures flow balance at non-root nodes, while \eqref{eq_topology6} forces $F_{ij}$ to zero when line $ij$ is not connected, utilizing the big-$M$ method.

\noindent (2) \textit{Number of PSES Constraint:}
\begin{flalign}
y_a^\mathrm{min} \leq \sum_{i\in \mathcal{I}_a} y_{i} \leq y_a^\mathrm{max}, \quad \forall a\in \mathcal{A}.
\label{eq_numberPSECS}
\end{flalign}
\vspace{-1em}

In cons. \eqref{eq_numberPSECS}, the total number of constructed PSESs is constrained to lie within the bounds of $y_a^\mathrm{min}$ and $y_a^\mathrm{max}$.

\noindent (3) \textit{Power Flow Constraints:}
\begin{flalign*}
\sum_{ij\in \mathcal{L}}P_{ij,t,s}-\sum_{ki\in \mathcal{L}} P_{ki,t,s}=P_{i,t}^\mathrm{Grid}(\xi)+y_i P_{i,t,s}^\mathrm{EV}-y_i P_{i,t,s}^\mathrm{PV}
\label{PowerFlow1}
\end{flalign*}
\vspace{-1.5em}
\begin{flalign}
&&+y_i\left(P_{i,t,s}^\mathrm{ESS,ch}-P_{i,t,s}^\mathrm{ESS,dch}\right), \quad \forall i\in \mathcal{I}, t\in \mathcal{T},s\in \mathcal{S},
\end{flalign}
\vspace{-1.5em}
\begin{flalign*}
\sum_{ij\in \mathcal{L}} Q_{ij,t,s}-\sum_{ki\in \mathcal{L}} Q_{ki,t,s}=Q_{i,t}^\mathrm{Grid}(\xi), 
\label{PowerFlow2}
\end{flalign*}
\vspace{-2em}
\begin{flalign}
&& \forall i\in \mathcal{I}, t\in \mathcal{T},s\in \mathcal{S},
\end{flalign}
\vspace{-2em}
\begin{flalign*}
{\left|V_{i,t,s}-V_{j,t,s}-2\left(R_{ij} P_{ij,t,s}+X_{ij} Q_{ij,t,s}\right)\right|} 
\label{PowerFlow3}
\end{flalign*}
\vspace{-2em}
\begin{flalign}
&&{\leq\left(1-z_{ij}\right) M, \quad \forall i,j\in \mathcal{I}, ij\in \mathcal{L}, t\in \mathcal{T}, s\in \mathcal{S}.}
\end{flalign}

The Linear DistFlow model in cons. \eqref{PowerFlow1}--\eqref{PowerFlow3} describes power flow in the radial DN. Cons. \eqref{PowerFlow1} ensures active power balance at node $i$, where the active power flowing out of node $i$ is balanced by the incoming power, considering the contributions from the grid-side load ($P_{i,t}^\mathrm{Grid}$), PV generation ($P_{i,t,s}^\mathrm{PV}$), EV charging demand ($P_{i,t,s}^\mathrm{EV}$), and the charging/discharging operations of the ESS ($P_{i,t,s}^\mathrm{ESS,ch}$, $P_{i,t,s}^\mathrm{ESS,dch}$). Similarly, cons. \eqref{PowerFlow2} enforces reactive power balance. {Cons. \eqref{PowerFlow3} captures the coupling relationship between squared voltage magnitudes and both active and reactive power flows on distribution lines \cite{wang2023joint}.}

\noindent (4) \textit{Grid Security Constraints:}
\begin{equation}
\left\lVert
\begin{array}{c}
P_{ij,t,s} \\
Q_{ij,t,s}
\end{array}
\right\rVert_2
\le z_{ij} S_{ij}^\mathrm{max},
\quad \forall\, ij \in \mathcal{L}, t \in \mathcal{T}, s \in \mathcal{S},
\label{eq_Security1}
\end{equation}
\vspace{-1em}
\begin{flalign}
V^\mathrm{min} \leq V_{i,t,s} \leq V^\mathrm{max}, \quad \forall i\in \mathcal{I}, t\in \mathcal{T}, s \in \mathcal{S}.
\label{eq_Security2}
\end{flalign}

Line capacity and voltage magnitude are limited by cons. \eqref{eq_Security1} and \eqref{eq_Security2} to ensure secure DN operation.

\noindent (5) \textit{Substation Power Constraints:}
\begin{flalign}
0 \leq P_{i,t,s}^\mathrm{Sub} \leq P^\mathrm{Sub,max}, \quad \forall i\in \mathcal{I}^\mathrm{Sub}, t \in \mathcal{T}, s \in \mathcal{S},
&\label{eq_Sub1}
\end{flalign}
\vspace{-2em}
\begin{flalign}
0 \leq Q_{i,t,s}^\mathrm{Sub}\leq Q^\mathrm{Sub,max}, \quad \forall i \in \mathcal{I}^\mathrm{Sub}, t \in \mathcal{T}, s \in \mathcal{S}.
&\label{eq_Sub2}
\end{flalign}

The active and reactive power flowing out of the substation are restricted by cons. \eqref{eq_Sub1} and cons. \eqref{eq_Sub2}, respectively.

\noindent (6) \textit{PV Operation Constraint:}
\begin{flalign}
0 \leq P_{i,t,s}^\mathrm{PV} \leq y_i \bar{P}_{i,t}^{\mathrm{PV}}(\xi), \quad \forall i\in \mathcal{I}, t\in \mathcal{T}, s \in \mathcal{S}.
&\label{eq_PV}
\end{flalign}

The output of PV is constrained by cons. \eqref{eq_PV} to remain within its maximum generation capacity.

\noindent (7) \textit{ESS Operation Constraints:}
\begin{flalign*}
{E_{i,t,s}^\mathrm{ESS}=\left(P_{i,t-1,s}^\mathrm{ESS,ch} \mu_i^\mathrm{ESS,ch}-P_{i,t-1,s}^\mathrm{ESS,dch} / \mu_i^\mathrm{ESS,dch}\right) \Delta t}
&\label{eq_ESS1}
\end{flalign*}
\vspace{-2em}
\begin{flalign}
&&{+E_{i,t-1,s}^\mathrm{ESS}, \quad \forall i \in \mathcal{I}, t \in \mathcal{T}, s \in \mathcal{S},}
\end{flalign}
\vspace{-2em}
\begin{flalign*}
P_i^\mathrm{ESS,min} \gamma_{i,t,s}^\mathrm{ESS,ch} \leq P_{i,t,s}^\mathrm{ESS,ch} \leq P_i^\mathrm{ESS,max} \gamma_{i,t,s}^\mathrm{ESS,ch}, 
&\label{eq_ESS2}
\end{flalign*}
\vspace{-2em}
\begin{flalign}
&& \forall i \in \mathcal{I}, t \in \mathcal{T}, s \in \mathcal{S},
\end{flalign}
\vspace{-2em}
\begin{flalign*}
P_i^\mathrm{ESS,min} \gamma_{i,t,s}^\mathrm{ESS,dch} \leq P_{i,t,s}^\mathrm{ESS,dch} \leq P_i^\mathrm{ESS,max} \gamma_{i,t,s}^\mathrm{ESS,dch},
&\label{eq_ESS3}
\end{flalign*}
\vspace{-2em}
\begin{flalign}
&& \forall i \in \mathcal{I}, t \in \mathcal{T}, s \in \mathcal{S},
\end{flalign}
\vspace{-2em}
\begin{flalign}
E_i^\mathrm{ESS,min} \leq E_{i,t,s}^\mathrm{ESS} \leq E_i^\mathrm{ESS,max}, \; \forall i \in \mathcal{I}, t \in \mathcal{T}, s \in \mathcal{S},
&\label{eq_ESS4}
\end{flalign}
\vspace{-2em}
\begin{flalign}
\gamma_{i,t,s}^\mathrm{ESS,ch}+\gamma_{i,t,s}^\mathrm{ESS,dch} \leq y_i, \quad \forall i \in \mathcal{I}, t \in \mathcal{T}, s \in \mathcal{S}.
&\label{eq_ESS5}
\end{flalign}

{Cons. \eqref{eq_ESS1}--\eqref{eq_ESS5} describe ESS operation. Each representative day comprises $T=24$ intervals of duration $\Delta t=1\,\mathrm{h}$, indexed by $\mathcal T=\{1,\ldots,T\}$. The variable $E_{i,t,s}^{\mathrm{ESS}}$ denotes stored energy at the beginning of interval $t$, and $E_{i,T+1,s}^{\mathrm{ESS}}$ denotes end-of-day energy. Cons. \eqref{eq_ESS1} updates stored energy using the charging and discharging powers in interval $t-1$, with respective efficiency coefficients $\mu_i^{\mathrm{ESS,ch}}$ and $\mu_i^{\mathrm{ESS,dch}}$. Cons. \eqref{eq_ESS2}--\eqref{eq_ESS4} bound the charging power, discharging power, and stored energy. Cons. \eqref{eq_ESS5} disables charging and discharging when no PSES is installed and prevents their simultaneous operation otherwise.
For each installed ESS, the initial condition $E_{i,1,s}^{\mathrm{ESS}}=E_i^{\mathrm{ESS},0}$ and the daily cyclic condition $E_{i,T+1,s}^{\mathrm{ESS}}=E_{i,1,s}^{\mathrm{ESS}}$ are imposed for all $s\in\mathcal S$, where $E_i^{\mathrm{ESS},0}$ is a prescribed initial energy within the admissible bounds. The cyclic condition prevents net depletion of stored energy over the representative day.}

\noindent (8) \textit{Carbon Emission Flow Model:}
\begin{flalign}
e_{i, t}=\frac{\displaystyle\sum_{ki \in \mathcal{L}} P_{ki,t} \rho_{ki,t}+P_{i,t}^\mathrm{G} e_{i,t}^\mathrm{G}}{\displaystyle\sum_{ki \in \mathcal{L}} P_{ki,t}+P_{i,t}^\mathrm{G}}, \; \forall i \in \mathcal{I} \backslash \mathcal{I}^\mathrm{Sub}, t \in \mathcal{T},
&\label{eq_CEF1}
\end{flalign}
\vspace{-1.2em}
\begin{flalign}
e_{i,t}^G=\frac{P_t^\mathrm{TG} e_t^\mathrm{TG}+P_t^\mathrm{TC} e_t^\mathrm{TC}}{P_t^\mathrm{TG}+P_t^\mathrm{TC}},\quad \forall i \in \mathcal{I}^\mathrm{Sub}, t \in \mathcal{T},
&\label{eq_CEF2}
\end{flalign}
\vspace{-1.5em}
\begin{flalign}
\boldsymbol{P}_{\mathrm{I}}\boldsymbol{E}_{\mathrm{I}}=\boldsymbol{P}_{\mathrm{L}}^{\top} \boldsymbol{E}_{\mathrm{I}}+\boldsymbol{P}_{\mathrm{G}}^{\top}\boldsymbol{E}_{\mathrm{G}},
&\label{eq_CEF3}
\end{flalign}
\vspace{-2em}
\begin{flalign}
\boldsymbol{P}_{\mathrm{I}}=\mathrm{diag}(
\boldsymbol{P}_{\mathrm{L}}^{\top} \mathbf{1}_{|I|}
+\boldsymbol{P}_{\mathrm{G}}^{\top} \mathbf{1}_{|G|}
),
&\label{eq_CEF4}
\end{flalign}
\vspace{-2em}
\begin{flalign*}
P_{i,t}^\mathrm{Load}=P_{i,t}^\mathrm{Grid}+P_{i,t}^\mathrm{EV}+P_{i, t}^\mathrm{ESS,ch}-P_\mathrm{i,t}^{PV}-P_{i,t}^\mathrm{ESS,dch},
\end{flalign*}
\vspace{-2em}
\begin{flalign}
&&\forall i \in \mathcal{I}, t \in \mathcal{T}.
\label{eq_CEF5}
\end{flalign}
\vspace{-1.5em}

Cons. \eqref{eq_CEF1}--\eqref{eq_CEF5} characterize the carbon emission flow model in the DN \cite{cheng2018planning}. Expression \eqref{eq_CEF1} defines the calculation of carbon intensity, where $\rho_{ki,t}$ denotes the carbon flow density on line $ki$, $P_{i,t}^\mathrm{G}$ is the active power injected into the system by generator sets, and $e_{i,t}^\mathrm{G}$ is the carbon emission intensity of the generator located at node $i$. The value of $e_{i,t}^\mathrm{G}$ is further determined by expression \eqref{eq_CEF2}, where $P_t^\mathrm{TG}$ and $P_t^\mathrm{TC}$ represent the power outputs of thermal and tidal current units, and $e_t^\mathrm{TG}$ and $e_t^\mathrm{TC}$ denote their respective carbon emission intensities.
By extending cons. \eqref{eq_CEF1}--\eqref{eq_CEF2} to the system-wide level, cons. \eqref{eq_CEF3} is obtained. In this expression, $\boldsymbol{P}_{\mathrm{I}}$, $\boldsymbol{P}_{\mathrm{L}}$, $\boldsymbol{P}_{\mathrm{G}}$, $\boldsymbol{E}_{\mathrm{I}}$, and $\boldsymbol{E}_{\mathrm{G}}$ denote the nodal active power flux matrix, branch power flow distribution matrix, power injection matrix, nodal carbon intensity vector, and generator carbon intensity vector, respectively. Expression \eqref{eq_CEF4} specifies the calculation of $\boldsymbol{P}_{\mathrm{I}}$, where $\mathbf{1}_{|I|}$ and $\mathbf{1}_{|G|}$ denote all-ones column vectors with dimensions $|I|$ and $|G|$, and $\mathrm{diag}(\cdot)$ forms the diagonal matrix. Finally, the active load at node $i$ is calculated by cons. \eqref{eq_CEF5}. {The CEF model is evaluated using the fixed network operating results for the retained planning design under the worst-case uncertainty identified by the upper problem, with the thermal--tidal supply mix determined by the lower-layer procurement problem. The power flow, PV/ESS dispatch, and EV demand quantities refer to the same evaluation condition, whose scenario dependence is suppressed for notational brevity. The matrices in \eqref{eq_CEF3}--\eqref{eq_CEF4} are constructed for each time interval.}

\noindent (9) \textit{Power Procurement Constraints:}
\begin{flalign}
P_t^\mathrm{TG}+P_t^\mathrm{TC}=\sum_{i \in \mathcal{I}} P_{i,t}^\mathrm{Load},
&\label{eq_EM1}
\end{flalign}
\vspace{-1.2em}
\begin{flalign}
0 \leq P_t^\mathrm{TG} \leq P_t^\mathrm{TG,max},
&\label{eq_EM2}
\end{flalign}
\vspace{-1.8em}
\begin{flalign}
0 \leq P_t^\mathrm{TC} \leq P_t^\mathrm{TC,max}.
&\label{eq_EM3}
\end{flalign}
\vspace{-1.5em}

Cons. \eqref{eq_EM1} ensures that the total power procured from thermal and tidal current units equals the active load. Cons. \eqref{eq_EM2} and \eqref{eq_EM3} impose their respective generation limits.

\vspace{-1em}
\subsection{Hybrid Uncertainty Modeling}
Uncertainty in this model arises from load fluctuations, PV availability, and EV charging behavior. EV charging uncertainty is represented by discrete spatio-temporal scenarios generated from charging profiles, while active load, reactive load, and PV availability are modeled as continuous uncertain quantities.

\noindent (1) \textit{Continuous Box Uncertainty for Load and PV Availability:}

The uncertainties of grid load and PV power output are represented using box-type uncertainty sets:
\begin{align}
&P_{i,t}^{\mathrm{Grid}}(\xi) \in \left[\bar{P}_{i,t}^{\mathrm{Grid}}-\Delta P_{i,t}^{-},\; \bar{P}_{i,t}^{\mathrm{Grid}}+\Delta P_{i,t}^{+}\right], \label{eq_boxP}\\
&Q_{i,t}^{\mathrm{Grid}}(\xi) \in \left[\bar{Q}_{i,t}^{\mathrm{Grid}}-\Delta Q_{i,t}^{-},\; \bar{Q}_{i,t}^{\mathrm{Grid}}+\Delta Q_{i,t}^{+}\right], \label{eq_boxQ}\\
&\bar{P}_{i,t}^{\mathrm{PV}}(\xi) \in \left[\bar{P}_{i,t}^{\mathrm{PV},0}-\Delta P_{i,t}^{\mathrm{PV},-},\; \bar{P}_{i,t}^{\mathrm{PV},0}+\Delta P_{i,t}^{\mathrm{PV},+}\right]. \label{eq_boxPV}
\end{align}

All continuous uncertainties are collected into \(\xi\in\Xi\), the Cartesian product of the above intervals.

\noindent (2) \textit{EV Charging Behaviors and Scenario Envelope:}

EVs can be incentivized through the TOU electricity tariff to adhere to scheduled charging profiles, thereby facilitating peak load shifting. For each scenario \(s\in\mathcal{S}\), given a set of EV sessions \(\mathcal{U}_s\) with a plug-in window \(\tau_u\subseteq\mathcal{T}\), the orderly charging profile can be obtained by minimizing the total charging payment:
\begin{flalign}
{\min_{\{P_{u,t,s}^{\mathrm{EV}}\}} \sum_{u \in \mathcal{U}_s} \sum_{t \in \tau_u} c^\mathrm{Ele}_t P_{u,t,s}^\mathrm{EV} \Delta t,}
&\label{eq_EV1}
\end{flalign}
\vspace{-1em}
\begin{flalign*}
\text{s.t.}\quad
{E_u^\mathrm{EV,min} \leq E_{u,s,0}^\mathrm{EV(\cdot)} + \sum_{t \in \tau_u} P_{u,t,s}^\mathrm{EV}\Delta t \leq E_u^\mathrm{EV,max},}
&\label{eq_EV2}
\end{flalign*}
\vspace{-2em}
\begin{flalign}
&& {\forall u \in \mathcal{U}_s,}
\end{flalign}
\vspace{-2em}
\begin{flalign}
{E_{u,s,0}^\mathrm{EV(\cdot)}+\sum_{t \in \tau_u} P_{u,t,s}^\mathrm{EV}\Delta t \geq E_{u,s}^\mathrm{EV(\cdot)}, \quad \forall u \in \mathcal{U}_s,}
&\label{eq_EV3}
\end{flalign}
\vspace{-1em}
\begin{flalign}
\sum_{u \in \mathcal{U}_s}\sum_{t \in \mathcal{T}\setminus\tau_u} P_{u,t,s}^\mathrm{EV}=0,
&\label{eq_EV4}
\end{flalign}
\vspace{-1em}
\begin{flalign}
0 \leq P_{u,t,s}^\mathrm{EV} \leq P_u^\mathrm{EV,max}, \quad \forall u \in \mathcal{U}_s, t \in \tau_u,
&\label{eq_EV5}
\end{flalign}
where $P_{u,t,s}^\mathrm{EV}$, $E_{u,s,0}^\mathrm{EV(\cdot)}$, and $E_{u,s}^\mathrm{EV(\cdot)}$ are the charging power of the $u$-th EV, its initial battery energy under uncertain behavior $(\cdot)$, and its target energy, respectively. Cons. \eqref{eq_EV2} limits the relationship between battery capacity and EV charging power. Cons. \eqref{eq_EV3} ensures that the target SoC is satisfied. Cons. \eqref{eq_EV4} enforces that EVs stop charging outside their plug-in windows, and the bounds of EV charging power are defined in cons. \eqref{eq_EV5}. 

Let \(i(u)\in\mathcal{I}\) denote the PSES candidate node to which session \(u\) is mapped. The scheduled EV demand aggregated at each node is
\begin{flalign}
P_{i,t,s}^{\mathrm{EV,sch}}
= \sum_{u\in\mathcal{U}_s:~i(u)=i} P_{u,t,s}^{\mathrm{EV}}, \quad \forall i,t,s.
\label{eq_EV6}
\end{flalign}
\vspace{-1em}

Meanwhile, to capture the possibility that EV users may not fully comply with the orderly charging schedule, a random aggregated profile $P_{i,t,s}^{\mathrm{EV,ran}}$ is constructed using a Monte Carlo method. The charging power of each EV session is modeled by a truncated normal distribution, and the following conservative envelope is adopted in implementation:
\begin{flalign}
P_{i,t,s}^{\mathrm{EV}} := \max\!\left\{P_{i,t,s}^{\mathrm{EV,ran}},~P_{i,t,s}^{\mathrm{EV,sch}}\right\},
\quad \forall i,t,s.
\label{eq_EV7}
\end{flalign}
\vspace{-1em}

The resulting $\{P_{i,t,s}^{\mathrm{EV}}\}_{i,t}$ defines the support point $\zeta_s^{\mathrm{EV}}$ for scenario $s$; $P_{i,t,s}^{\mathrm{EV}}$ serves as the scenario input in \eqref{PowerFlow1}.

\noindent (3) \textit{Discrete EV Scenarios and Ambiguity Set for Probabilities:}

EV charging uncertainty is represented by a finite scenario set \(\{\zeta^{\mathrm{EV}}_s\}_{s\in\mathcal{S}}\), where each scenario is a spatio-temporal charging profile
\(\zeta^{\mathrm{EV}}_s := \{P_{i,t,s}^{\mathrm{EV}}\}_{i\in\mathcal{I},t\in\mathcal{T}}\).
The scenario probabilities \(\pi := (\pi_s)_{s\in\mathcal{S}}\) are not assumed known exactly; an ambiguity set centered at a reference estimate \(\pi^0\) \cite{zhao2015data} is adopted instead:
\begin{align}
\Pi^{\mathrm{EV}}
:= \Big\{ \pi \in \mathbb{R}_{+}^{|\mathcal{S}|}\; \Big| \;
&\sum_{s\in\mathcal{S}}\pi_s=1, \label{eq:pi1}\\
&\sum_{s\in\mathcal{S}}|\pi_s-\pi_s^0|\le \theta_1, \label{eq:pi2}\\
&|\pi_s-\pi_s^0|\le \theta_\infty,\;\forall s\in\mathcal{S}
\Big\}, \label{eq:pi3}
\end{align}
where the radii \(\theta_1,\theta_\infty\) are set based on confidence levels \(\alpha_1,\alpha_\infty\) and EV sample size \(N^{\mathrm{EV}}\) (Appendix~\ref{app_1} and \ref{app_2}).

\vspace{-1em}
\section{Solution Method with i-C\&CG Algorithm}

\subsection{Compact Formulation of the Optimization Model}
Based on the preceding physical and uncertainty models, the planning problem can be structured as a tri-layer model. The upper layer determines DN topology and PSES siting, the mid-layer optimizes network operation under uncertainty, and the lower layer performs carbon accounting and power procurement in response to the coastal city's low-carbon policies and renewable-resource conditions.
For clarity, the compact form of the proposed model is summarized below and detailed in \eqref{eq_hol_1}--\eqref{eq_hol_6}.
\begin{flalign}
\min _{\mathbf{z},\mathbf{y}} \left\{ f(\mathbf{z},\mathbf{y}) + \max _{\boldsymbol{\delta}} \min _{\mathbf{p},\mathbf{q}} g(\mathbf{p},\mathbf{q}) + \min _{\mathbf{e},\mathbf{p}^{\prime}} h(\mathbf{e},\mathbf{p}^{\prime})\right\},
\label{eq_hol_1}
\end{flalign}
\vspace{-2em}
\begin{flalign}
\text{s.t.}\quad \mathbf{A}\mathbf{z}+\mathbf{B}\mathbf{y} \geq \mathbf{b},
\label{eq_hol_2}
\end{flalign}
\vspace{-2em}
\begin{flalign}
\mathbf{C}\mathbf{z}+\mathbf{D}\mathbf{y}+\mathbf{d} \geq \mathbf{E}\mathbf{p}+\mathbf{F}\mathbf{q}+\mathbf{G}\boldsymbol{\delta},
\label{eq_hol_3}
\end{flalign}
\vspace{-2em}
\begin{flalign}
\mathbf{H}\mathbf{z}+\mathbf{I}\mathbf{y}+\mathbf{r} \geq \mathbf{J}\mathbf{p}+\mathbf{J}^{\prime}\mathbf{p}^{\prime}+\mathbf{K}\mathbf{q}+\mathbf{L}\mathbf{e}+\mathbf{M}\boldsymbol{\delta}
\label{eq_hol_4}
\end{flalign}
\vspace{-2em}
\begin{flalign}
\mathbf{z} \in \{0,1\}^{n_z},\mathbf{y} \in \{0,1\}^{n_y},\boldsymbol{\delta}:=(\xi,\pi)\in \Xi\times\Pi^{EV},
\label{eq_hol_5}
\end{flalign}
\vspace{-2em}
\begin{flalign}
\mathbf{p} \in \mathbb{R}^{n_{p}},\mathbf{p}^{\prime} \in \mathbb{R}^{n_{p^{\prime}}},\mathbf{q} \in \mathbb{R}^{n_{q}},\mathbf{e} \in \mathbb{R}_{+}^{n_{e}},
\label{eq_hol_6}
\end{flalign}
where $\mathbf{z}$ and $\mathbf{y}$ represent the decision variables associated with DN lines and PSESs; $\mathbf{p}$ and $\mathbf{p}^{\prime}$ stack the mid-/lower-layer active power variables, and $\mathbf{q}$ stacks the reactive power variables. $\mathbf{e}$ is the system-wide carbon intensity vector, and $\boldsymbol{\delta}$ captures uncertainties. The upper-layer objective $f(\mathbf{z},\mathbf{y})$ in \eqref{eq_hol_1} is formulated in equation \eqref{eq_Inv} and constrained by \eqref{eq_hol_2}--\eqref{eq_hol_3}, which correspond to cons. \eqref{eq_topology1}--\eqref{eq_numberPSECS} and \eqref{PowerFlow1}--\eqref{eq_ESS5}. The mid-layer operation objective $g(\mathbf{p},\mathbf{q})$ aligns with equation \eqref{eq_Ope}, sharing constraints in \eqref{eq_hol_3}. At the lower layer, $h(\mathbf{e},\mathbf{p}^{\prime})$ is defined by equations \eqref{eq_Carb}--\eqref{eq_Powr}, addressing carbon emissions and power procurement costs, subject to \eqref{eq_hol_4} corresponding to cons. \eqref{eq_CEF1}--\eqref{eq_EM3}. $\mathbf{A}$--$\mathbf{M}$ denote coefficient matrices, while $\mathbf{b}$, $\mathbf{d}$, and $\mathbf{r}$ are the corresponding vectors in the constraints. 

 \begin{figure}[htbp]
    \centering
    \includegraphics[width=0.40\textwidth]{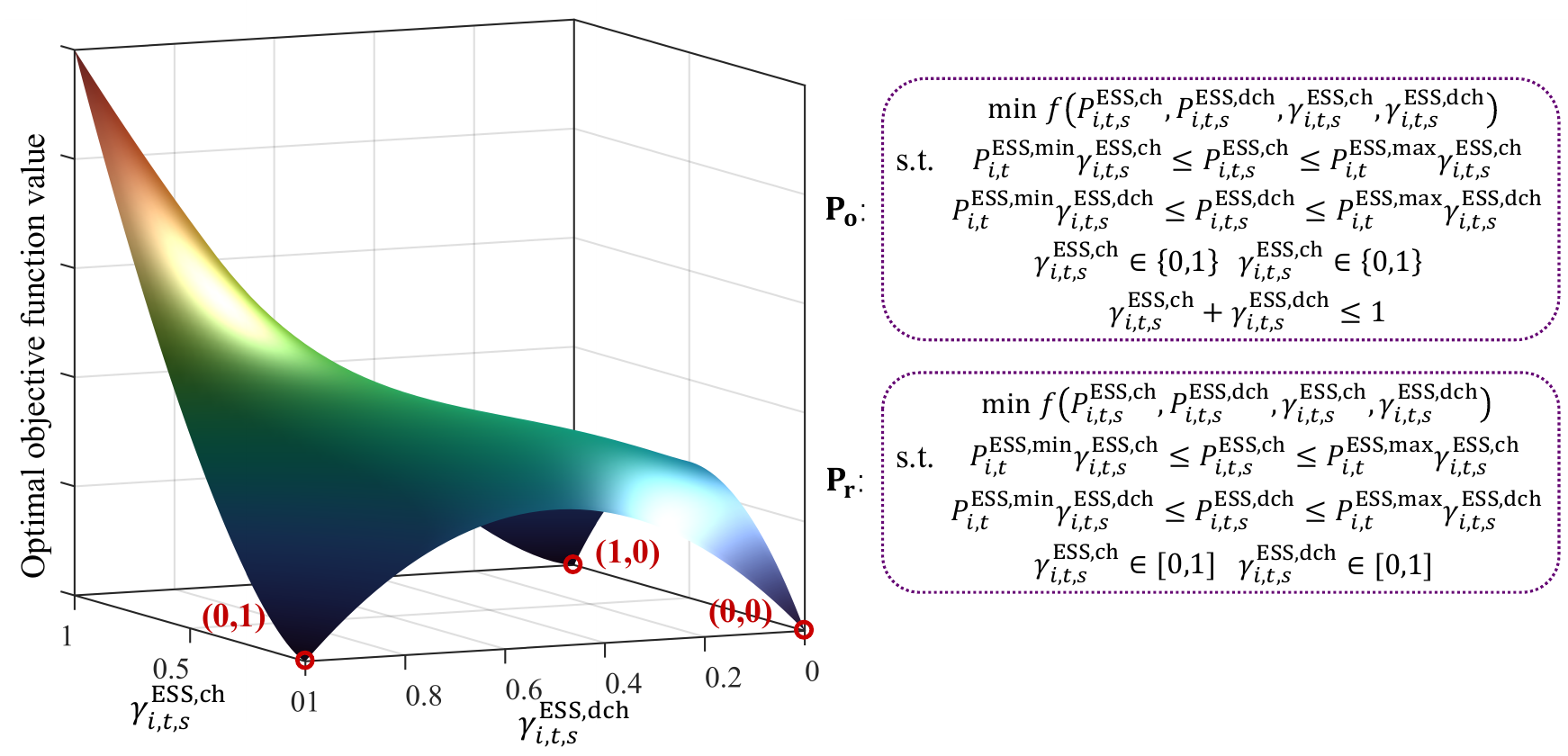}
    \vspace{-1.0em}
    \caption{A small example for feasibility analysis.} \label{fig_ESS}
\end{figure}

Due to the presence of binary decision variables in the ESS-related constraints \eqref{eq_ESS1}--\eqref{eq_ESS5} in the mid-layer, the problem is strongly non-convex and computationally intensive. As a result, it is not suitable for a direct KKT-based reformulation, which requires a continuous convex problem. To resolve this challenge, we relax the binary variables $(\gamma_{i,t,s}^{\mathrm{ESS,ch}},\,\gamma_{i,t,s}^{\mathrm{ESS,dch}})$ from $\{0,1\}^2$ to $[0,1]^2$ and prove the relaxation is exact, following the exact-relaxation rationale in \cite{li2015sufficient}.  

\noindent \textbf{Proposition 1.} Fix any PSES siting decision $y$ and assume the minimum admissible ESS charging and discharging powers are zero. Let $\mathbf{P_o}$ denote the original mid-layer problem with the ESS mode indicators constrained as $(\gamma_{i,t,s}^{\mathrm{ESS,ch}},\,\gamma_{i,t,s}^{\mathrm{ESS,dch}})\in\{0,1\}^2$, and let $\mathbf{P_r}$ denote its relaxed counterpart with $(\gamma_{i,t,s}^{\mathrm{ESS,ch}},\,\gamma_{i,t,s}^{\mathrm{ESS,dch}})\in[0,1]^2$. Then the relaxation is exact in the sense that $\mathbf{P_o}$ and $\mathbf{P_r}$ attain the same optimal objective value. Moreover, $\mathbf{P_r}$ admits an optimal solution whose ESS mode indicators are binary.

\noindent \textbf{Remark 1.} Unlike \cite{li2015sufficient}, which focuses on complementarity constraints imposed on continuous charging/discharging power variables in the operational stage with pre-installed ESS (without siting binaries), our formulation includes the unavoidable PSES siting variable $y_i$, and thus explicitly gates the ESS mode indicators via \eqref{eq_ESS5}, deactivating ESS operations when $y_i=0$.

\noindent \textbf{Proof.} Consider $\mathbf{P_o}$ and its relaxation $\mathbf{P_r}$ in Fig.~\ref{fig_ESS}, let $\mathcal{F}_o$ and $\mathcal{F}_r$ denote their feasible sets.
The only difference between the two problems lies in the domain of the mode indicators in \eqref{eq_ESS5}; all other ESS constraints \eqref{eq_ESS1}--\eqref{eq_ESS4} remain identical.
Therefore, $\mathcal{F}_o \subseteq \mathcal{F}_r$.

Let $v_o^\star := \min_{x\in\mathcal{F}_o} f(x),
v_r^\star := \min_{x\in\mathcal{F}_r} f(x)$,
where $x$ collects the decision variables of the mid-layer problem, including the ESS variables.
Since $\mathcal{F}_o \subseteq \mathcal{F}_r$, it immediately follows that $v_r^\star \le v_o^\star$.

Next, we show that $\mathbf{P_r}$ has an optimal solution with binary mode indicators. Let $x_r^\star$ be an optimal solution of $\mathbf{P_r}$ with objective value $v_r^\star$. When $y_i=0$, constraint \eqref{eq_ESS5} forces $\gamma_{i,t,s}^{\mathrm{ESS,ch}}=\gamma_{i,t,s}^{\mathrm{ESS,dch}}=0$,
so the corresponding ESS is inactive and the indicators are already binary.
Now consider any $(i,t,s)$ with $y_i=1$. Because the ESS is lossy with $\mu_i^{\mathrm{ESS,ch}}<1$ and $\mu_i^{\mathrm{ESS,dch}}>1$, simultaneous charging and discharging within the same time interval represents a purely dissipative behavior. It introduces additional round-trip losses without providing any operational advantage for minimizing the mid-layer objective. Hence, an optimal solution to the relaxed problem can always be selected so that charging and discharging do not occur simultaneously at each $(i,t,s)$. Under such a solution, at least one of the two mode indicators equals zero for every $(i,t,s)$.
Given this property and the zero lower charging/discharging-power bound, we can recover binary mode indicators without changing feasibility or the objective value: for each $(i,t,s)$, set $\gamma_{i,t,s}^{\mathrm{ESS,ch}}=1$ if $P_{i,t,s}^{\mathrm{ESS,ch}}>0$ and $0$ otherwise, and set $\gamma_{i,t,s}^{\mathrm{ESS,dch}}=1$ if $P_{i,t,s}^{\mathrm{ESS,dch}}>0$ and $0$ otherwise. Then \eqref{eq_ESS5} remains satisfied, and the power bounds \eqref{eq_ESS2}--\eqref{eq_ESS3} continue to hold. This yields a feasible solution of $\mathbf{P_o}$ that attains the same objective value $v_r^\star$. Consequently, $v_o^\star \le v_r^\star$.

Combining $v_r^\star \le v_o^\star$ and $v_o^\star \le v_r^\star$ gives $v_o^\star=v_r^\star$.
Therefore, the relaxation is exact and the ESS retains its three operating modes, namely charging, discharging, and idle.
With Proposition~1, the mid-layer becomes a continuous convex program, so its KKT conditions can be employed to reformulate the original bi-level ``$\max \min$'' problem into a single-level ``$\max$'' problem, as shown in \eqref{eq_kkt_1}--\eqref{eq_kkt_7}.
\vspace{-0.7em}
\begin{flalign}
\max _{{\boldsymbol{\delta}},\mathbf{p},\mathbf{q},\boldsymbol{\lambda},\boldsymbol{\omega}} \hspace{0.5em}
g(\mathbf{p},\mathbf{q}),
\label{eq_kkt_1}
\end{flalign}
\vspace{-2.0em}
\begin{flalign}
\text{s.t.}\quad \mathbf{C}\mathbf{z}+\mathbf{D}\mathbf{y}-\mathbf{E}\mathbf{p}-\mathbf{F}\mathbf{q}-\mathbf{G}\boldsymbol{\delta}+\mathbf{d} \geq 0,
\label{eq_kkt_2}
\end{flalign}
\vspace{-2.0em}
\begin{flalign}
\mathbf{C}\mathbf{z}+\mathbf{D}\mathbf{y}-\mathbf{E}\mathbf{p}-\mathbf{F}\mathbf{q}-\mathbf{G}\boldsymbol{\delta}+\mathbf{d} \leq M(\mathbf{1}-\boldsymbol{\omega}),  
\label{eq_kkt_3}
\end{flalign}
\vspace{-2.2em}
\begin{flalign}
0 \leq \boldsymbol{\lambda}\leq M\boldsymbol{\omega},
\label{eq_kkt_4}
\end{flalign}
\vspace{-2.2em}
\begin{flalign}
\nabla_{\mathbf{p}} \hspace{0.2em} g(\mathbf{p},\mathbf{q}) + \mathbf{E}^\mathbf{T}\boldsymbol{\lambda} = 0,
\label{eq_kkt_5}
\end{flalign}
\vspace{-2.2em}
\begin{flalign}
\nabla_\mathbf{q} \hspace{0.2em} g(\mathbf{p},\mathbf{q}) + \mathbf{F}^\mathbf{T}\boldsymbol{\lambda} = 0, 
\label{eq_kkt_6}
\end{flalign}
\vspace{-2.2em}
\begin{flalign}
\boldsymbol{\lambda} \in \mathbb{R}^m_{+}, \boldsymbol{\omega} \in \{0,1\}^m,
\label{eq_kkt_7}
\end{flalign}
where $\boldsymbol{\lambda}$ denotes the Lagrange multiplier, $\boldsymbol{\omega}$ denotes the auxiliary binary variables introduced for complementary slackness, and $m$ is the number of inequality constraints in the KKT conditions. Cons. \eqref{eq_kkt_2} are the primal constraints, while cons. \eqref{eq_kkt_3}--\eqref{eq_kkt_6} contain the big-$M$ linearized complementary slackness conditions and the dual constraints \cite{zeng2013solving}.

\noindent \textbf{Proposition 2.} The third-layer problem for carbon emission and power procurement is downstream of the planning and operation layers and can be decoupled from them.
{After the planning decisions and the network operating inputs used for evaluation have been fixed, the third-layer variables are optimized without feeding back into the upper two layers.}

\noindent \textbf{Remark 2.}
The CEF is a virtual flow evaluated from the realized power-flow solution and the associated supply mix, capturing how carbon emissions are allocated across generation, transmission, and consumption. As CEF is derived based on the results of system planning and operation, it does not enter or modify the physical feasibility region of the DN in the upper two layers \cite{cheng2018planning},\cite{yang2023improved}.
Similarly, the coastal city procurement decision only allocates the realized net demand among available sources and does not feed back into the planning or operational constraints.
Therefore, once the upper two layers yield the solutions $\mathbf{z}^*$, $\mathbf{y}^*$, $\mathbf{p}^*$, $\mathbf{q}^*$, and $\boldsymbol{\delta}^*$, they can be substituted into \eqref{eq_holiCE_2}, leading to the decoupled third-layer formulation \eqref{eq_holiCE_1}--\eqref{eq_holiCE_2} with decision variables $\mathbf{p}^{\prime}$ and $\mathbf{e}$ only:
\vspace{-0.5em}
\begin{flalign}
\min_{\mathbf{e},\mathbf{p}^{\prime}} \hspace{0.5em}h(\mathbf{e},\mathbf{p}^{\prime}),
\label{eq_holiCE_1}
\end{flalign}
\vspace{-1.9em}
\begin{flalign}
\text{s.t.}\quad \mathbf{H}\mathbf{z^*}+\mathbf{I}\mathbf{y^*}+\mathbf{r}\geq \mathbf{J}\mathbf{p^*}+\mathbf{J}^{\prime}\mathbf{p}^{\prime}+\mathbf{K}\mathbf{q^*}+\mathbf{L}\mathbf{e}+\mathbf{M}\boldsymbol{\delta^*}.
\label{eq_holiCE_2}
\end{flalign}

\vspace{-2.0em}
\subsection{Inexact Column-and-Constraint Generation Algorithm}

After decoupling the third layer, the upper two layers are reformulated as a two-stage DRO problem. The upper layer passes the DN topology $\mathbf{z}^*$ and PSES siting decisions $\mathbf{y}^*$ to the mid-layer. In return, the mid-layer identifies the worst-case scenario ${\mathbf{\delta}^*}$ and the corresponding scenario-dependent constraints \eqref{eq_kkt_2}--\eqref{eq_kkt_6}. {This reformulated problem is then solved using the i-C\&CG algorithm, with problem-specific master and operational subproblem formulations for topology-and-siting co-planning.}
The algorithm is initialized in Step 1 by setting the necessary parameters. In each iteration, the master problem is solved with a relative optimality gap $\varepsilon_{i}^\mathrm{MP}$, yielding the upper bound $UB_i$ and the lower bound $LB_i$.
If the solution is valid ($LB_i \geq \overline{LB}$), the index $\ell$ is updated to the current iteration to indicate the most recently valid solution. $\overline{LB}$ is then updated to $UB_i$ to accelerate the improvement of the lower bound in the next iteration of the master problem.
With the optimal solution from the master problem fixed, the subproblem is solved to identify the worst-case scenario $\mathbf{\delta}_i^*$ and its corresponding optimal value, following the traditional C\&CG procedure. If $\overline{UB}$ and $UB_i$ are sufficiently close, the algorithm enters the exploitation phase with a backtracking routine (Step 11), which leverages the most recent valid lower bound $LB_\ell$ and returns to Step 3 to re-solve the master problem based on $LB_\ell$ using the reduced relative gap $\alpha\varepsilon_{i}^\mathrm{MP}$ to correct the inaccuracies in the previous iteration.
Otherwise, the algorithm proceeds to the exploration phase by enlarging the scenario set as in the C\&CG algorithm. The algorithm terminates when the convergence criterion ${(\overline{UB}-LB_\ell)}/{\overline{UB}} < \varepsilon$ is satisfied.
\begin{figure*}[htbp]
    \centering
    \includegraphics[width=0.90\textwidth]{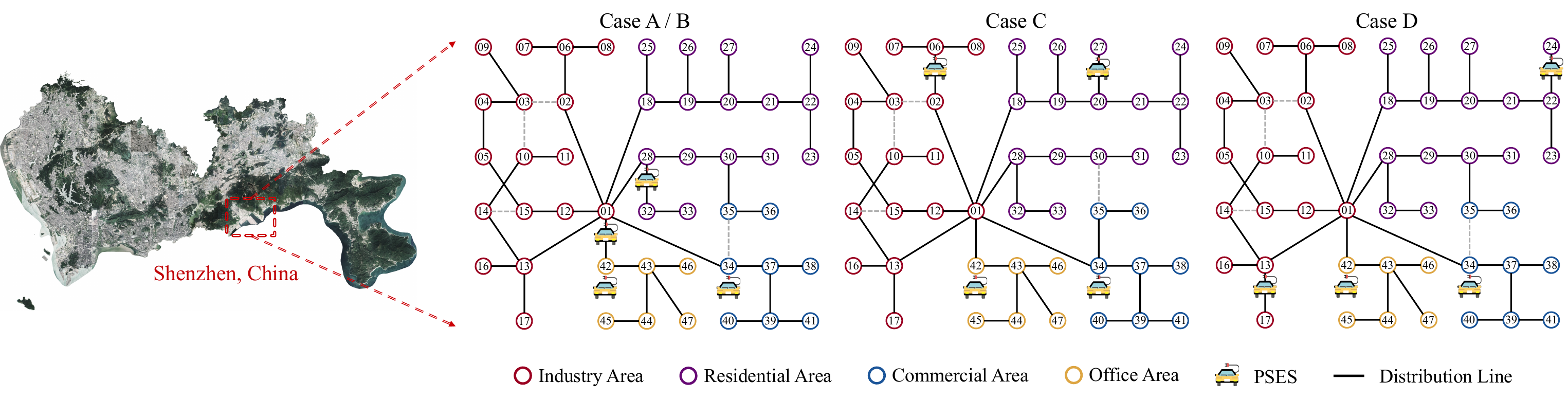}
    \vspace{-1.0em}
    \caption{Planning solutions for different cases in Shenzhen, China.}\label{fig_SZ}
\end{figure*}

\begin{table}[htbp]
\footnotesize
\centering
\renewcommand{\arraystretch}{1.1}
\begin{tabular}{@{}>{\raggedleft\arraybackslash}p{2em}@{\hspace{0.5em}}p{0.92\linewidth}@{}}
\Xhline{1pt}
\multicolumn{2}{@{}l}{\textbf{Algorithm 1:}  \quad i-C\&CG Algorithm for Solving the Proposed Problem} \\
\hline
1: & \textbf{Initialize:} $\overline{LB} = -\infty$, $\overline{UB} = +\infty$, $i = 1$, $\ell = 0$, $\mathcal{D} = \emptyset$, $\varepsilon \in [0,1]$, \\
& $\tilde{\varepsilon} \in \left(0, {\varepsilon}/{(\varepsilon + 1)}\right)$ , $\varepsilon_{i}^\mathrm{MP}$, $\alpha \in (0,1)$. \\
2: & \textbf{while} ${(\overline{UB}-LB_\ell)}/{\overline{UB}} \geq \varepsilon$  \textbf{do}\\
3: & \quad Within a relative optimality gap $\varepsilon_{i}^\mathrm{MP}$, solve the \textbf{Master Problem}:\\
& \qquad \qquad \qquad $v_i^*= \min_{\mathbf{z},\mathbf{y}} f(\mathbf{z},\mathbf{y}) + \eta$ \\
& \qquad \qquad \qquad \qquad \quad $\text{s.t.} \quad \eta \geq g(\mathbf{p},\mathbf{q})$, \\
& \qquad \qquad \qquad \qquad \qquad \quad $f(\mathbf{z},\mathbf{y}) + \eta \geq \overline{LB}$, \\
& \qquad \qquad \qquad \qquad \qquad \quad Cons. \eqref{eq_hol_2}, \eqref{eq_kkt_2}--\eqref{eq_kkt_6}. \\
4: & \quad Record the optimal solution $(\mathbf{z}_i^*,\mathbf{y}_i^*, \eta_i^*)$. \\
5: & \quad Record a valid lower bound $LB_i \geq \overline{LB}$  and an upper bound\\
& \quad $UB_i=v_i^*$. If $LB_i \geq \overline{LB}$, then set $\ell \gets i$.\\
6: & \quad Update $\overline{LB} =UB_i$. \\
7: & \quad For fixed $\mathbf{z} = \mathbf{z}_i^*$ and $\mathbf{y} = \mathbf{y}_i^*$, solve the \textbf{Subproblem}: \\
& \qquad \qquad \qquad $D_i^*= \max_{\mathbf{\delta},\mathbf{p},\mathbf{q},\boldsymbol{\lambda},\boldsymbol{\omega}} g(\mathbf{p},\mathbf{q})$ \\
& \qquad \qquad \qquad \qquad \quad $\text{s.t.} \quad$ Cons. \eqref{eq_kkt_2}--\eqref{eq_kkt_6}. \\
8: & \quad Record the worst scenario $\mathbf{\delta}_i^*$ and its optimal value $D_i^*$.\\
9: & \quad Update $\overline{UB} = \min\{\overline{UB}, f(\mathbf{z}_i^*,\mathbf{y}_i^*)+ D_i^*\}$. \\
10: & \quad \textbf{if} ${(\overline{UB} - UB_i)}/{\overline{UB}} < \tilde{\varepsilon}$  \textbf{then} \\
11: & \quad \quad \textbf{Exploitation:} Set $i \gets \ell$, $\overline{LB} = LB_\ell$, update $\varepsilon_{i}^\mathrm{MP} = \alpha  \varepsilon_{i}^\mathrm{MP}$\\
& \quad \quad for all $i \geq \ell$ with \textbf{Backtracking Routine}.  \\
12: & \quad \quad \textbf{break} \\
13: & \quad \textbf{else} \\
14: & \quad \quad \textbf{Exploration:} Enlarge the scenario set $\mathcal{D} = \mathcal{D} \cup \{\mathbf{\delta}_i^*\}$, update\\
& \quad \quad $i \gets i + 1$. \\
15: & \quad \textbf{end if} \\
16: & \textbf{end while} \\
17: & \textbf{Output:} Return the optimal solution ($\mathbf{z}_{i-1}^*$,$\mathbf{y}_{i-1}^*$) for which \\
&$obj^*=\overline{UB}$.\\
\hline
\end{tabular}
\end{table}

\vspace{-1.0em}
\section{Case Studies}

{The proposed planning framework is tested on a 47-node planning network derived from the real-world DN case in Shenzhen, China, whose detailed network parameters are publicly available in \cite{gao2026longgangdata}. }
{The 47-node network is divided into four functional areas, $\mathcal{A}=\{1,2,3,4\}$: industrial (Area~1, Nodes~1--17), residential (Area~2, Nodes~18--33), commercial (Area~3, Nodes~34--41), and office (Area~4, Nodes~42--47).}
The investment costs of DN lines and conventional EV charging stations (EVCSs) are 23.30 $\times 10^4$ CNY/km and 187 $\times 10^4$ CNY/unit, while that for the PSES is 344.5 $\times 10^4$ CNY/unit, comprising a 75 kW PV system and an ESS with 1.5 MWh storage, 0.2 MW charging, and 0.3 MW discharging power. PV installations receive a 1 CNY/W upfront subsidy and 0.05 CNY/kWh energy subsidy, while ESSs receive 0.9 CNY/Wh installation and 12 CNY/kWh capacity standby subsidies. 
Representative coastal conditions are assumed: the annual-mean relative humidity and ambient temperature are set to 80\% and $22^{\circ}\mathrm{C}$, and the chloride deposition rate is 300 mg/(m\textsuperscript{2}$\cdot$ day). The protective thicknesses $d^{\mathrm{Zn}}$ and $d^{\mathrm{Al}}$ are set to 45 $\mu$m and 25 $\mu$m \cite{ISO9223-2012}.
{Substituting these values into the dose--response expressions, gives estimates of approximately $\kappa_{\mathrm{Line}}^{\mathrm{Salt}}=0.13$ and $\kappa_{\mathrm{PSES}}^{\mathrm{Salt}}=0.09$ for DN lines and PSESs, respectively. Under the assumed spatially homogeneous coastal conditions, each coefficient is applied uniformly within its respective asset category.}
All components use a 0.05 discount rate, with design lives of 20 years for DN lines and EV charging stations, 30 years for PV, and 25 years for ESS. The voltage magnitude is constrained between 0.9 and 1.1 p.u.
The thermal unit has a 30 MW capacity, 0.85 t/MWh carbon intensity, and 0.4 CNY/kWh procurement cost. Installed tidal capacity is 9 MW; other parameters are detailed in Fig. \ref{fig_Tidal}(a). Shenzhen's TOU tariffs are 1.11 CNY/kWh (peak), 0.65 CNY/kWh (off-peak), and 0.25 CNY/kWh (standard) \cite{wang2026integrated}. Carbon prices are 50/40/30 CNY/t for industrial zones, commercial/office areas, and residential areas. For EV scenario generation, $\alpha_1$ and $\alpha_\infty$ are set to 0.99, standard deviation $\sigma$ to 0.12, and the confidence level $1-\alpha$ to 0.95.

\vspace{-1.0em}
\subsection{Simulation Results}
To isolate the effects of PSES deployment, salt-spray modeling, and the i-C\&CG algorithm, four cases are presented.
The model was formulated using the YALMIP tool in MATLAB (2021a) and evaluated with GUROBI Optimizer (9.5.2) on the Apple M3 Pro (12-core CPU, 18-core GPU). The cases are outlined as follows:

{\textit{Case A}: Coastal DN planning with conventional EVCSs considering salt-spray corrosion via the C\&CG algorithm}

{\textit{Case B}: Coastal DN planning with conventional EVCSs considering salt-spray corrosion via the i-C\&CG algorithm}

{\textit{Case C}: Coastal DN--PSES planning without considering salt-spray corrosion via the i-C\&CG algorithm}

{\textit{Case D}: Coastal DN--PSES planning considering salt-spray corrosion via the i-C\&CG algorithm (proposed model)}

The convergence gap is set to 1\%, {and one charging station is required to be deployed in each area and four stations in total.}
{The four cases support three pairwise comparisons. Cases A and B solve the same salt-spray-aware DN–EVCS planning problem using the C\&CG and i-C\&CG algorithms, respectively, thereby isolating the algorithmic effect at the identical convergence tolerance. Cases B and D both consider salt-spray corrosion and use the i-C\&CG algorithm; they compare conventional EVCS planning with integrated PSES planning under their respective equipment, investment, and subsidy assumptions. Cases C and D both use PSESs and the i-C\&CG algorithm and differ in whether salt-spray degradation costs are included, allowing the influence of this modeling choice on the computed planning solutions to be assessed. Network topology and station locations remain optimization outputs in all cases. Thus, B–D evaluates the complete station-technology alternative, whereas C–D evaluates the inclusion of the salt-spray cost component.
{For each case, the carbon intensities, the carbon emission and power procurement costs are calculated using the network operating results associated with the worst-case
uncertainty identified for the retained planning design.}
The planning results of DN topology and station siting, along with the corresponding costs and solution time, are presented in Fig. \ref{fig_SZ} and Table \ref{tab:comparison results}.}

\begin{table}[htbp]
  
  \centering
  \caption{Annualized Costs($10^{4}$ CNY) and Solution Time(s)}
  \vspace{-1em}
  \label{tab:comparison results}
  \begin{tabular}{ccccc}
  \hline\hline
  Case No. & \textit{Case A} & \textit{Case B} & \textit{Case C} & \textit{Case D*}  \\ \hline
  DN      &$\surd$ &$\surd$ &$\surd$ &$\surd$   \\ 
  Conventional EVCS    &$\surd$ &$\surd$ & --- & ---  \\
  PSES    & --- & --- &$\surd$ &$\surd$  \\
  Salt-spray corrosion  &$\surd$ &$\surd$ & --- &  $\surd$   \\
  \hline
  i-C\&CG algorithm     &  ---   &$\surd$   &$\surd$   &$\surd$   \\ 
  C\&CG algorithm  &$\surd$ & --- & --- & --- \\ 
 \hline
  Line construction cost &59.58 &59.58 &59.76 &59.58   \\ 
  Station investment cost  &60.00 &60.00 &104.00 &104.00     \\ 
  Salt-spray degradation cost  &13.15 &13.15   &0.00  &17.11     \\
  Network loss cost    &44.43 &44.43  &41.89  &43.38     \\ 
  Subsidy  &0.00 &0.00 &-51.60  &-51.60  \\  \hline
  Cost of upper two layers   &177.16 &177.16 &154.05  &172.47    \\ \hline
  Carbon emission cost  &225.46 &225.46 &218.63  &217.08  \\ 
  Power procurement cost  &3256.88 &3256.88 &3237.90  &3237.89  \\ 
 \hline
  Total cost    &3659.50 &3659.50 &3610.58  &3627.44    \\ \hline
  Solution time  &5772.65 &585.87 &9310.34  &8599.13  \\ 
  \hline\hline
  \end{tabular}
\end{table}

In the planning results, since Case A and Case B differ only in the solution algorithm, they yield identical DN topologies and EVCS siting, with selected nodes located near the substation to reduce power losses. In Case C, the DN topology is reconfigured (Lines 34--35 and 30--35), and four PSESs are split between the substation vicinity and the feeder end. This result shows that the loss reduction from PV and ESS must be balanced against the additional losses induced by EV charging. In Case D, three PSESs are sited close to the substation and one at the feeder end (Node 24), indicating that the loss reduction from PV and ESS outweighs the incremental losses associated with EV charging. A comparison of the four cases shows that modeling salt-spray corrosion and deploying PV/ESS can substantially alter the optimal DN topology and PSES siting decisions.
Comparing the annualized results, Case D achieves a 2.36\% reduction in power loss cost compared to Cases A and B, attributed to the effect of PV and ESS in PSESs. Meanwhile, carbon emission and power procurement costs decrease by 3.72\% and 0.58\%, respectively, resulting in a total annualized cost reduction of 32.06$\times 10^4$ CNY. Compared with Case C ignoring salt-spray corrosion, the salt-spray degradation cost in Case D accounts for 9.92\% of the upper two layers' cost, raising the total cost by $16.86\times 10^4$ CNY.
The line construction cost and network loss cost also differ across the cases, further indicating that coastal DN planning requires a more refined representation of coastal stressors and their impacts.
Under the same convergence tolerance, Cases A and B share identical cost outcomes. {However, the solution time of Case B using the i-C\&CG algorithm is only 10.15\% of that of Case A with the conventional C\&CG algorithm, illustrating the computational benefit of the adopted i-C\&CG implementation in this benchmark comparison.}

\begin{figure}[htbp]
    \centering
    \includegraphics[width=0.45\textwidth]{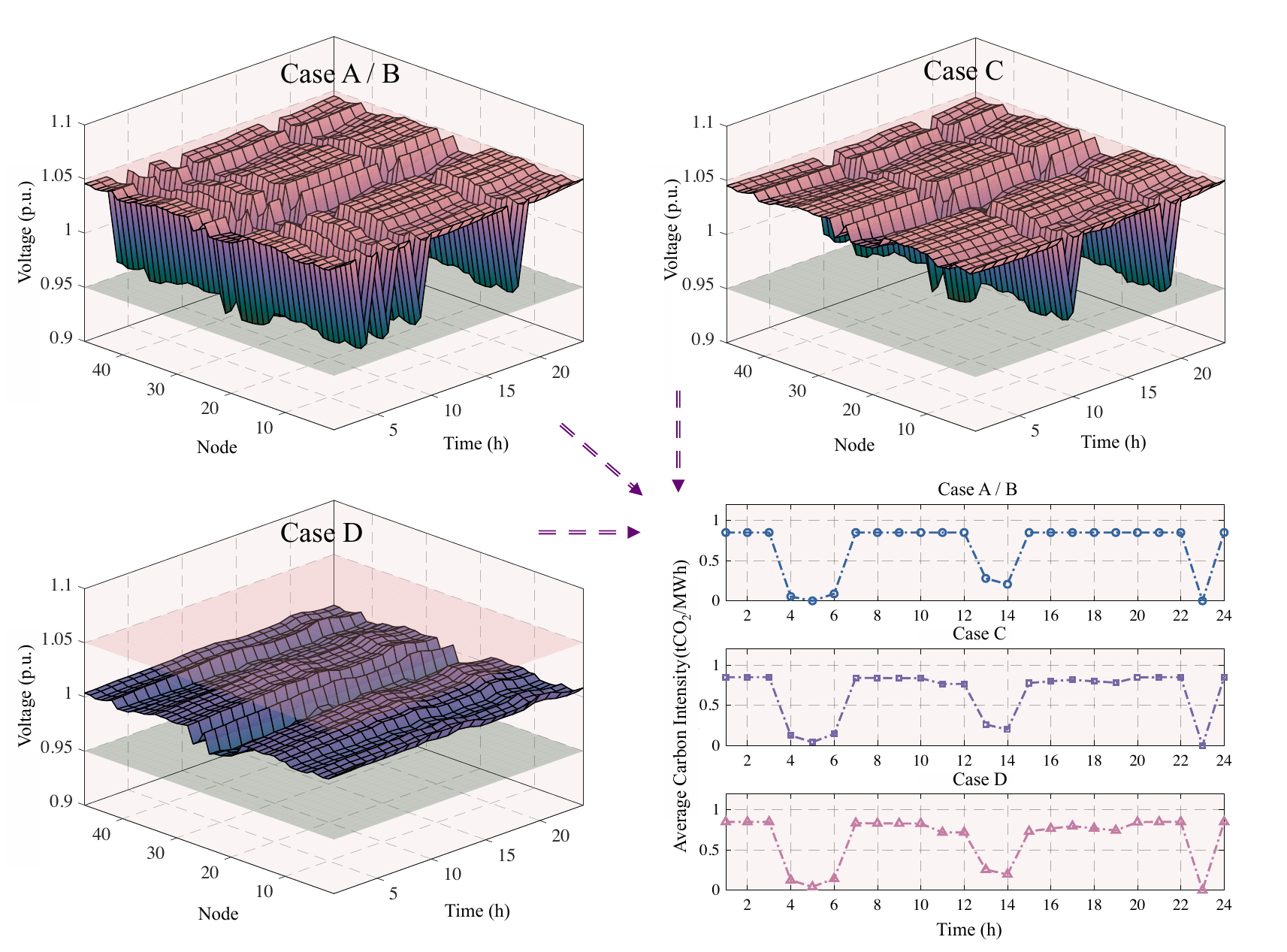} 
    \vspace{-1.0em}
    \caption{{Voltage and node-averaged carbon intensity profiles for different cases.}}\label{fig_vol}
\end{figure}

The nodal voltage and average carbon intensity for different cases are shown in Fig.~\ref{fig_vol}. The voltage ranges for Case A/B, Case C, and Case D are 0.1000, 0.1000, and 0.0315 p.u., with standard deviations of 0.0298, 0.0222, and 0.0064 p.u., respectively. Case D exhibits the smallest voltage range and deviation, indicating the most stable voltage profile and the highest voltage quality. Compared with conventional EV charging stations, PSES-integrated PV and ESS generally smooth voltage fluctuations.
The average carbon intensity across all nodes and time intervals for Cases A/B, C, and D is 0.6636, 0.6494, and 0.6367 tCO\textsubscript{2}/MWh, which indicates that integrating PV and ESS within PSESs can effectively reduce nodal carbon intensity. Compared to siting PSESs near the substation, constructing them at feeder-end nodes will lead to lower carbon intensity of the surrounding downstream nodes, but this configuration will increase power losses. Therefore, planning should consider the trade-off among multiple cost components to achieve an optimal balance.

\begin{figure}[htbp]
    \centering
    \includegraphics[width=0.49\textwidth]{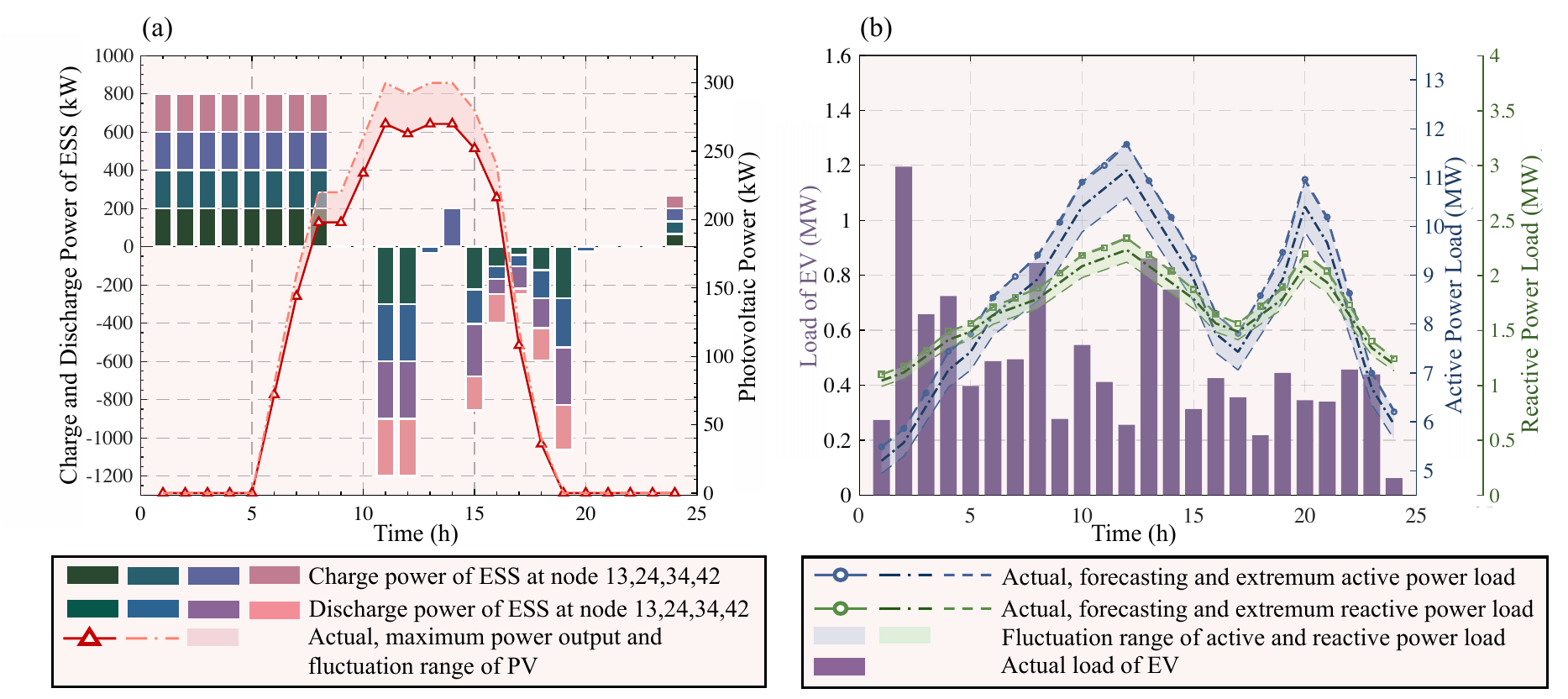}
    \vspace{-1.5em}
    \caption{(a) Operation of PV and ESS. (b) Load fluctuation errors.}\label{fig_PQEV}
\end{figure}

The operational characteristics of PV and ESS, along with load fluctuation errors in Case D, are shown in Fig. \ref{fig_PQEV}. The active and reactive power loads exhibit similar fluctuation trends throughout the day. The ESS discharges during peak load intervals and charges during low load intervals, demonstrating its effectiveness in peak shaving and valley filling within the DN. The PV availability reaches its lower uncertainty bound of 0.9 and the available PV power is fully utilized in all time intervals, indicating that the proposed model remains feasible under the worst-case scenario. The EV charging load exhibits fluctuations similar to those of the active and reactive loads, highlighting the effectiveness of electricity price--driven orderly charging strategies for EV users.

\vspace{-1em}
\subsection{Sensitivity Analysis of Different Parameters}

\begin{figure}[htbp]
    \centering
    \includegraphics[width=0.45\textwidth]{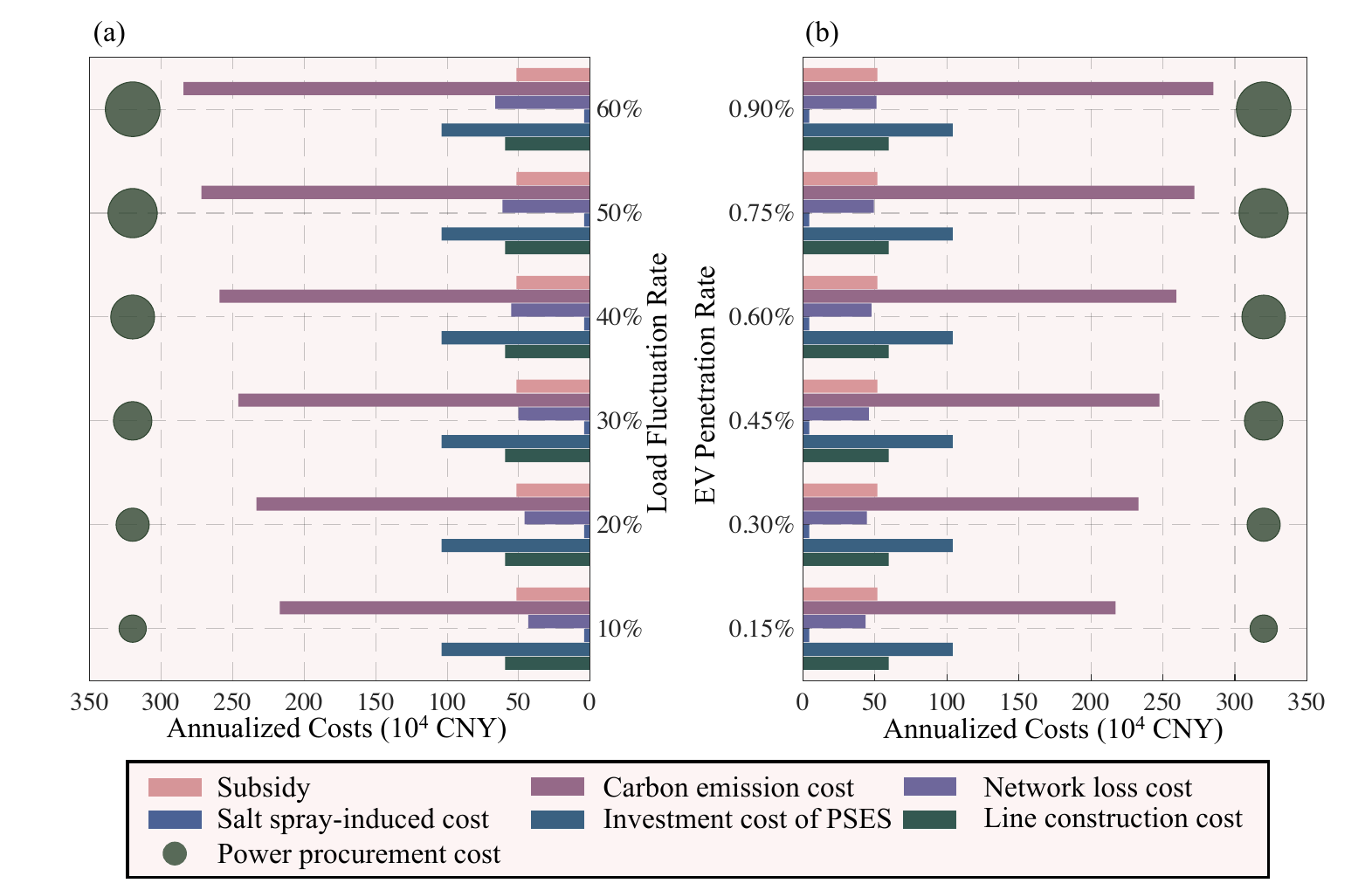}
    \vspace{-0.7em}
    \caption{Sensitivity analysis of load fluctuation and EV penetration.}\label{fig_Sensitivity}
\end{figure}

Fig. \ref{fig_Sensitivity} shows the variations in different cost components with increasing load fluctuation and EV penetration defined in \cite{wang2023joint}. As both load and EV penetration increase, power losses and total costs rise accordingly. However, the DN topology remains unchanged, confirming the planning framework's robustness to future growth in EV penetration.

\begin{figure}[htbp]
    \centering
    \includegraphics[width=0.45\textwidth]{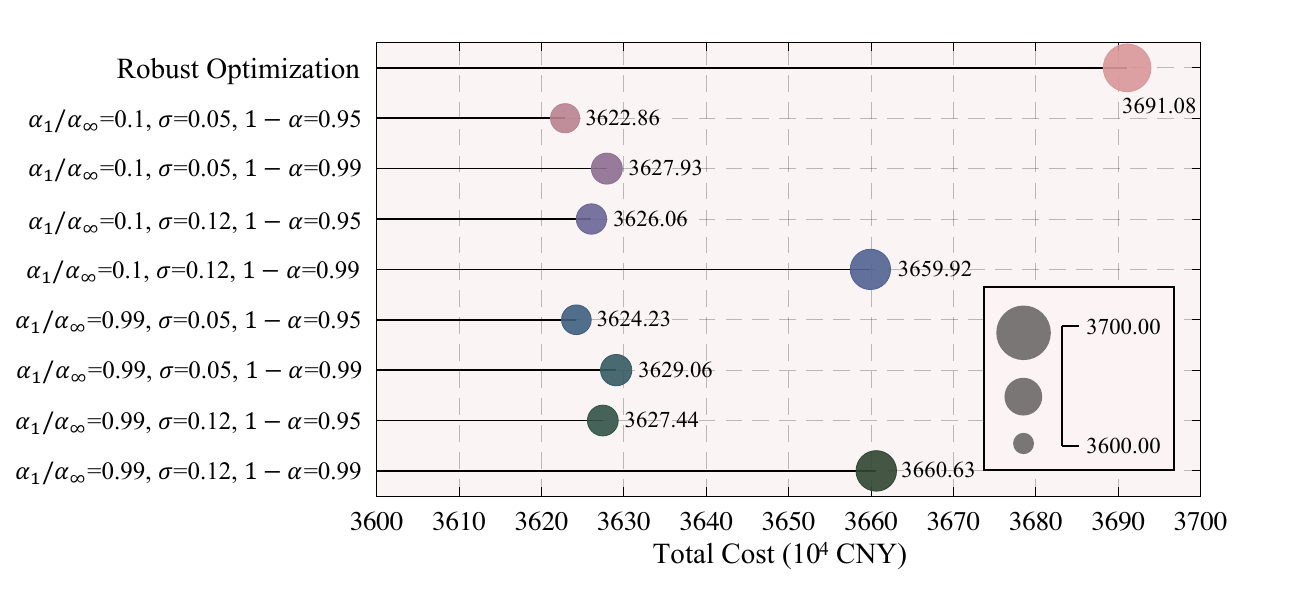}
    \vspace{-1em}
    \caption{Sensitivity analysis of distributionally robust parameters.}\label{fig_DRO}
\end{figure}

Fig.~\ref{fig_DRO} reports the sensitivity of key RO and DRO parameters. Since RO optimizes for the worst case, it yields the most conservative solution and the highest total cost. In DRO, larger $\alpha_1$ and $\alpha_\infty$ indicate increased ambiguity in EV scenario probabilities. Similarly, higher values of the standard deviation $\sigma$ and confidence level $1-\alpha$ used to model EV random charging power in the truncated normal distribution imply increased variability. Since DRO considers the worst-case distribution within an ambiguity set, the total cost escalates as the related uncertainty parameters grow. Compared to RO, DRO more accurately reflects the empirical distribution of EV charging data while preserving robustness against adverse scenarios.

\vspace{-1em}

\vspace{-0.5em}
\subsection{Analysis of Tidal Current Units in Coastal Cities}

\begin{figure}[htbp]
    \centering
    \includegraphics[width=0.49\textwidth]{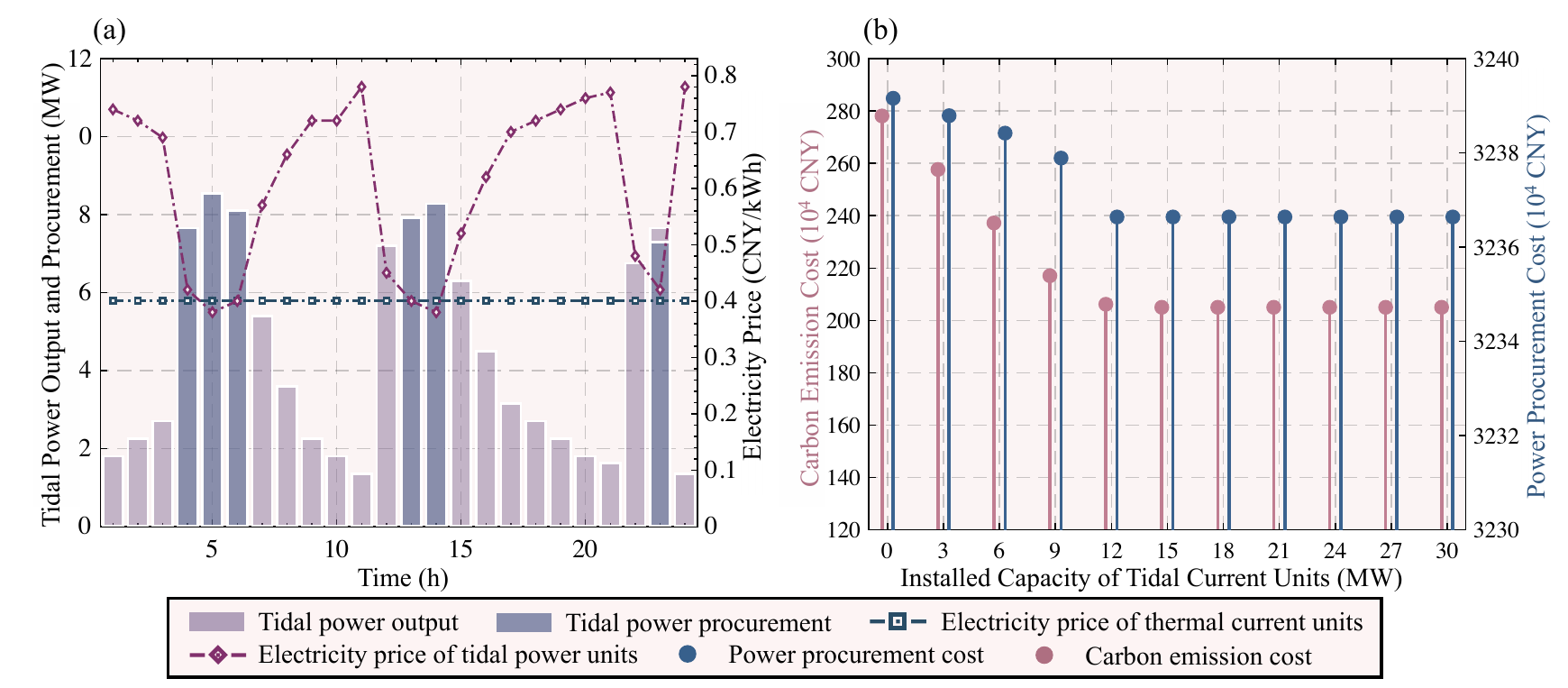}
    \vspace{-2em}
    \caption{(a) Power output, procurement and electricity price. (b) Sensitivity analysis of the installed capacity of tidal current units.}\label{fig_Tidal}
\end{figure}

Tidal energy output varies with tidal velocity, peaking between high and low tides under the semi-diurnal tide model. Although higher tidal output generally lowers its electricity price, the levelized cost of tidal units remains above that of thermal generators due to their high life-cycle costs. As shown in Fig. \ref{fig_Tidal}(a), tidal power is preferred for procurement when its electricity price is relatively low. Conversely, when the tidal electricity price significantly exceeds that of thermal power, it will not be selected. Fig. \ref{fig_Tidal}(b) presents the sensitivity analysis of installed tidal capacity. Increasing capacity up to 15 MW markedly reduces both carbon emissions and procurement costs, lowering the total system cost. However, once the installed capacity exceeds 15 MW, the tidal pricing mechanism prevents further expansion from increasing tidal power procurement because it no longer contributes to cost minimization. These results indicate that the economics of marine renewable procurement should be considered in coastal DN--PSES planning for low-carbon energy transition.

\vspace{-1.0em}

\section{Conclusion}

This paper has proposed a data-driven distributionally robust planning framework for coastal DNs and PSESs in support of the low-carbon energy transition. The framework jointly optimizes radial DN topology and PSES siting while incorporating salt-spray degradation costs, tidal-generation-driven procurement costs, and CEF-based nodal carbon accounting. It therefore extends conventional DN planning by linking coastal physical stressors, transport electrification, distributed PV-storage flexibility, and low-carbon operation in one planning model.

To address heterogeneous uncertainty, the model combines box uncertainty for continuous load and PV deviations with probability ambiguity over EV-charging scenarios generated from historical charging profiles. The resulting two-stage DRO problem is reformulated by relaxing the ESS mode binaries without loss of optimality and applying a KKT-based decomposition. A tailored i-C\&CG algorithm is then used to reduce the computational burden caused by coupled topology and siting binaries. 

{Case studies on a realistic 47-node coastal feeder show that salt-spray exposure can change both DN topology and PSES siting decisions. Compared with conventional EV charging station planning, the proposed DN--PSES co-planning improves voltage quality and reduces the annualized total cost by \(32.06\times 10^4\) CNY in the benchmark case. The sensitivity analysis further shows that tidal procurement changes the cost--carbon trade-off and that its economic benefit becomes limited beyond moderate installed capacity. Compared with conventional C\&CG, the tailored i-C\&CG algorithm obtains the same planning design and cost within the specified convergence tolerance, with about an order-of-magnitude reduction in solution time.} 

{These results indicate that coastal DN--PSES planning for low-carbon energy transition should jointly consider charging infrastructure, network topology, distributed PV-storage flexibility, marine renewable procurement, carbon accounting, and coastal degradation. In particular, coastal distribution planners should assess salt-spray degradation costs together with local PV–storage flexibility when comparing network and charging-station configurations. These factors can alter the investment–operation trade-off and thus affect which configuration is economically preferred under the adopted planning assumptions. 
The current salt-spray model uses annual-average environmental conditions and spatially uniform degradation coefficients within each asset category, which limits its representation of location-specific exposure.
Future work will extend the framework to broader coastal renewable portfolios, such as offshore wind, and to spatially varying coastal environmental conditions.}
\vspace{-1em}
\appendices

\section{EV Scenario Generation Method}
\label{app_1}

We generate discrete EV scenarios \(\{\zeta_s^{\mathrm{EV}}\}_{s\in\mathcal{S}}\) from historical charging data and construct \(\pi_s^0\) from scenario frequencies.

\noindent \subsubsection{Data and Preprocessing}
We assume a dataset of EV charging sessions. Each session is mapped to a PSES candidate bus \(i(u)\) and discretized into a time-slot trajectory over \(\mathcal{T}\).
For each day \(d=1,\ldots,N^{\mathrm{EV}}\), we compute an \emph{scheduling} aggregate profile by solving \eqref{eq_EV1}--\eqref{eq_EV6} (with \(s\leftarrow d\)) and a \emph{random} aggregate profile via Monte Carlo, then form the conservative envelope \eqref{eq_EV7} to obtain \(\mathbf{p}^{(d)}:=\{p^{(d)}_{i,t}\}_{i,t}\).

\noindent \subsubsection{Clustering and Scenario Construction}
We cluster daily profiles \(\{\mathbf{p}^{(d)}\}_{d=1}^{N^{\mathrm{EV}}}\) into \(|\mathcal{S}|\) clusters using K-means. Each cluster centroid yields a representative scenario profile \(\zeta_s^{\mathrm{EV}}\).

\noindent \subsubsection{Reference Probability \(\pi^0\)}
Let \(N_s\) be the number of days assigned to cluster \(s\). The reference probability is
\begin{align}
\pi_s^0 = \frac{N_s}{N^{\mathrm{EV}}}, \quad \forall s\in\mathcal{S}.
\label{eq:pi0}
\end{align}

The detailed pseudocode is shown in Algorithm~2.
\begin{table}[!htbp]
\footnotesize
\centering
\renewcommand{\arraystretch}{1.1}
\begin{tabular}{@{}>{\raggedleft\arraybackslash}p{2em}@{\hspace{0.5em}}p{0.92\linewidth}@{}}
\Xhline{1pt}
\multicolumn{2}{@{}l}{\textbf{Algorithm 2:} \quad EV Scenario Generation Method} \\
\hline
1: & \textbf{Input:} Historical charging sessions, \(\mathcal{T}\), \(N^{\mathrm{EV}}\), \(|\mathcal{S}|\). \\
2: & \textbf{for} $d=1,\ldots,N^{\mathrm{EV}}$ \textbf{do} \\
3: & \quad Build session set $\mathcal{U}^{(d)}$, windows $\tau_u$, and energy targets. \\
4: & \quad Solve the scheduling charging LP \eqref{eq_EV1}--\eqref{eq_EV6} to obtain $P^{\mathrm{EV,sch}}_{i,t}(d)$. \\
5: & \quad Construct the random profile $P^{\mathrm{EV,ran}}_{i,t}(d)$ with Monte Carlo method. \\
6: & \quad Set envelope $p^{(d)}_{i,t}=\max\{P^{\mathrm{EV,sch}}_{i,t}(d),P^{\mathrm{EV,ran}}_{i,t}(d)\}$. \\
7: & \textbf{end for} \\
8: & Cluster profiles $\{\mathbf{p}^{(d)}\}$ into $|\mathcal{S}|$ to obtain centroids $\zeta^{\mathrm{EV}}_s$. \\
9: & Set $\pi_s^0 = N_s/N^{\mathrm{EV}}$ for scenario $s$. \\
10: & \textbf{Output:} scenario set $\{\zeta_s^{\mathrm{EV}}\}$, reference probabilities $\pi_s^0$. \\
\hline
\end{tabular}
\end{table}
\vspace{-1em}

\section{Ambiguity Radii and Confidence Levels}
\label{app_2}
Following data-driven DRO concentration bounds, we set:
\begin{align}
\theta_1=\frac{|\mathcal{S}|}{2N^{\mathrm{EV}}}\ln\frac{2|\mathcal{S}|}{1-\alpha_1}, \qquad
\theta_\infty=\frac{1}{2N^{\mathrm{EV}}}\ln\frac{2|\mathcal{S}|}{1-\alpha_\infty}.
\label{eq:theta}
\end{align}

These expressions link \(\theta_1\) and \(\theta_\infty\) to the sample size \(N^{\mathrm{EV}}\) and confidence levels \(\alpha_1,\alpha_\infty\), while the combined \(1\)- and \(\infty\)-norm constraints bound deviations of \(\pi\) from \(\pi^0\) \cite{zhao2015data}.

\bibliographystyle{IEEEtran}

\small\bibliography{clean_reference}

\end{document}